\documentclass[aop,seceqn,citesort,MSNbibl,dvips]{arximspdf}
\usepackage{subenv}

\doi{10.1214/09-AOP501}
\volume{38}
\issue{3}
\pubyear{2010}
\firstpage{1180}
\lastpage{1220}

\makeatletter

\newcommand{\mazinti}[1]{\mbox{\fontsize{8.36}{9.36}\selectfont{#1}}}
\newcommand{\mmazinti}[1]{\mbox{\fontsize{6.6}{7.6}\selectfont{#1}}}

\newproclaim{convention}{Convention}[section]

\newtheorem{theorem}[convention]{Theorem}

\newproclaim{remark}[convention]{Remark}
\newproclaim{conjecture}[convention]{Conjecture}

\newtheorem{lemma}[convention]{Lemma}
\newtheorem{corollary}[convention]{Corollary}
\newtheorem{proposition}[convention]{Proposition}

\newcommand{\bDelta}{\bolds\Delta}

\makeatother

\begin{document}
\begin{frontmatter}

\title{Optimal local H\"{o}lder index for density states of
superprocesses with $(1+\beta)$-branching mechanism\protect\thanksref{T1}}
\runtitle{Optimal local H\"{o}lder index for superprocess states}


\begin{aug}
\author[A]{\fnms{Klaus} \snm{Fleischmann}\corref{}\ead[label=e1]{fleischm@wias-berlin.de}},
\author[B]{\fnms{Leonid} \snm{Mytnik}\ead[label=e2]{leonid@ie.technion.ac.il}\ead[url,label=u1]{http://ie.technion.ac.il/leonid.phtml}} and
\author[C]{\fnms{Vitali} \snm{Wachtel}\ead[label=e3]{wachtel@mathematik.uni-muenchen.de}}
\runauthor{K. Fleischmann, L. Mytnik and V. Wachtel}
\affiliation{Weierstrass Institute, Technion Israel Institute of
Technology and~University~of~Munich}
\address[A]{K. Fleischmann\\
Weierstrass Institute\\
\quad for Applied Analysis\\
\quad and Stochastics\\
Mohrenstrasse 39\\
D-10117 Berlin\\
Germany\\
\printead{e1}} 
\address[B]{L. Mytnik\\
Faculty of Industrial Engineering\\
\quad and Management\\
Technion Israel Institute\\
\quad of Technology\\
Haifa 32000\\
Israel\\
\printead{e2}\\
\printead{u1}}
\address[C]{V. Wachtel\\
Mathematical Institute\\
University of Munich\\
Theresienstrasse 39\\
D-80333 Munich\\
Germany\\
\printead{e3}}
\end{aug}

\thankstext{T1}{Supported by the German Israeli Foundation for
Scientific Research and
Development, Grant G-807-227.6/2003.}

\received{\smonth{5} \syear{2008}}
\revised{\smonth{4} \syear{2009}}

%
\begin{abstract}
For $ 0<\alpha\leq2$, a super-$\alpha$-stable motion $ X$ in
$ \mathsf{R}^{d}$ with branching of index $ 1+\beta\in(1,2)$ is
considered. Fix arbitrary $t>0$. If $ d<\alpha/\beta$, a dichotomy for
the density function of the measure $ X_{t}$ holds: the density
function is locally H\"{o}lder continuous if $ d=1$ and
$ \alpha>1+\beta$ but locally unbounded otherwise. Moreover, in the
case of continuity, we determine the optimal local H\"{o}lder index.
\end{abstract}

%

\end{frontmatter}

\section{Introduction and statement of results}

\subsection{Background and purpose}\label{SS.background}

For $ 0<\alpha\leq2$, a super-$\alpha$-stable motion
$X=\{X_{t}\dvtx t\geq0\}$ in $\mathsf{R}^{d}$ with branching of index
$ 1+\beta\in(1,2]$ is a finite measure-valued process related to
the log-Laplace equation%
%
%
\begin{equation} \label{logLap}%
\frac{d}{dt}\,u =
\bDelta
_{\alpha}u +au- bu^{1+\beta},
\end{equation}
where $ a\in\mathsf{R}$ and $ b>0$ are any fixed constants. Its
underlying motion is described by the fractional Laplacian
$\bDelta_{\alpha}:=-(- \bDelta)^{\alpha/2}$ determining a symmetric
$\alpha$-stable motion in $\mathsf{R}^{d}$ of index $ \alpha\in(0,2]$
(Brownian motion if $ \alpha=2)$ whereas its continuous-state branching
mechanism described by
%
%
\begin{equation} \label{not.Psi}%
v \mapsto-av+bv^{1+\beta} =: \Psi(v),\qquad v\geq0,
\end{equation}
belongs to the domain of attraction of a stable law of index $ 1+\beta
\in(1,2]$ (the branching is critical if $ a=0$).

It is well known that in dimensions $ d<\frac{\alpha}{\beta}$ at
any fixed time $ t>0$ the measure $ X_{t}=X_{t}(dx)$
is absolutely continuous with probability one (cf. Fleischmann
\cite{Fleischmann1988critical} where $a=0;$ the noncritical case requires
the obvious changes). By an abuse of notation, we sometimes denote a version
of the density function of the measure $X_{t}=X_{t}(dx)$ by the same
symbol, $ X_{t}(dx)=X_{t}(x)\, dx$, that is,
$ X_{t}=\{X_{t}(x)\dvtx x\in\mathsf{R}^{d}\}$. In the case
of one-dimensional continuous super-Brownian motion ($\alpha=2$, $\beta=1$),
even a joint-continuous density field $ \{X_{t}(x)\dvtx t>0, x\in
\mathsf{R}\}$ exists, satisfying a stochastic equation (Konno
and Shiga \cite{KonnoShiga1988} as well as Reimers \cite{Reimers1989}).

From now on we assume that $d<\frac{\alpha}{\beta}$ and
$\beta\in(0,1)$. For the Brownian case $ \alpha=2$ and
if $a=0$ (critical branching), Mytnik \cite{Mytnik2002} proved that a version
of the density $\{X_{t}(x)\dvtx t>0, x\in\mathsf{R}^{d}\}$ of the
measure $X_{t}(dx)\,dt$ exists that satisfies, in a weak sense,
the following stochastic partial differential equation (SPDE):
%
%
\begin{equation}
\frac{\partial}{\partial t}X_{t}(x)=%
\bDelta X_{t}(x)+(bX_{t-}(x))^{1/(1+\beta)}\dot{L}(t,x),
\end{equation}
where $\dot{L}$ is a $(1+\beta)$-stable noise without negative jumps.
\begin{convention}
\label{Conv}
From now on, (if it is not stated otherwise explicitly) we
use the term \textit{density} to denote the density function of the
measure $X_{t}(dx)$ with respect to the Lebesgue measure.
\end{convention}

For the same model (as in the paragraph before Convention \ref{Conv}), in
Mytnik and Perkins \cite{MytnikPerkins2003} regularity and irregularity
properties of the density at fixed times had been revealed. More precisely,
these densities have continuous versions if $ d=1$, whereas they
are locally unbounded on open sets of positive $X_{t}(dx)$-measure in
all higher dimensions $ (d<\frac{2}{\beta} )$.

The first \textit{purpose} in the present paper is to allow also discontinuous
underlying motions, that is to consider also all $ \alpha\in(0,2)$.
Then actually the same type of \textit{fixed time dichotomy} holds
(recall that $ d<\frac{\alpha}{\beta})$: continuity of densities
if $ d=1$ and $ \alpha>1+\beta$ whereas local
unboundedness is true if $ d>1$ or $ \alpha\leq1+\beta$.

However, the \textit{main purpose} of the paper is to address the following
question: what is the optimal local H\"{o}lder index in the first case of
existence of a continuous density? Here by optimality we mean that
there is a
critical index $\eta_{\mathrm{c}}$ such that for any fixed $t>0$ there
is a
version of the density which is locally H\"{o}lder continuous of any index
$\eta<\eta_{\mathrm{c} }$ whereas there is no locally H\"{o}lder continuous
version with index $ \eta\geq\eta_{\mathrm{c} }$.

In \cite{MytnikPerkins2003} continuity of the density at fixed times is proved
by some moment methods, although moments of order larger than $1+\beta$
are in
general infinite in the $1+\beta<2$ case. A standard procedure to get local
H\"{o}lder continuity is the Kolmogorov criterion by using
``high'' moments. This, for instance, can be done in the
$\beta=1$ case ($\alpha=2$, $d=1)$ to show local H\"{o}lder continuity
of any
index smaller than $\frac{1}{2}$ (see the estimates in the proof of
Corollary 3.4 in Walsh \cite{Walsh1986}).

Due to the lack of ``high'' moments in our $\beta<1$ case we cannot use
moments to get the optimal local H\"{o}lder index. Therefore we have to
get deeply into the jump structure of the superprocess to obtain the
needed estimates. As a result we are able to show the \textit{local
H\"{o}lder continuity} of all orders $ \eta<\eta
_{\mathrm{c}}:=\frac{\alpha}{1+\beta}-1$, provided that $ d=1$ and $
\alpha>1+\beta$. We also verify that the bound $ \eta_{\mathrm{c}}$ for
the local H\"{o}lder index is in fact \textit{optimal} in the sense
that there are points $x_{1},x_{2}$ such that
the density increments $ |X_{t}(x_{1})-X_{t}(x_{2})|$
are of a larger order than $ |x_{1}-x_{2}|^{\eta}$ as
$ x_{1}-x_{2}\rightarrow0$ for every $ \eta\geq\eta
_{\mathrm{c} }$. For precise formulations, see
Theorem \ref{T.prop.dens} below.

\subsection{Statement of results}\label{SS.statement}

Write $ \mathcal{M}_{\mathrm{f}}$ for the set of all finite measures
$\mu$ defined on $ \mathsf{R}^{d}$ and $ |\mu |$ for its total mass
$\mu(\mathsf{R}^{d})$. Let $ \Vert f\Vert_{U}$ denote the essential
supremum (with respect to Lebesgue measure) of a function $
f\dvtx\mathsf{R}^{d}\rightarrow\mathsf
{R}%
_{+}:=[0,\infty)$ over a nonempty open set $ U\subseteq
\mathsf{R}^{d}$.

Let $ p^{\alpha}$ denote the continuous $\alpha$\textit{-stable
transition kernel} related to the fractional Laplacian $ %
\bDelta_{\alpha}=-(-\bDelta)^{\alpha/2}$, and $ S^{\alpha}$ the related
\textit{semigroup}.

Recall that $ 0<\alpha\leq2$, $ 1+\beta\in(1,2)$ and
$ d<\frac{\alpha}{\beta}$, and consider again the
$(\alpha,d,\beta)$-superprocess $ X=\{X_{t}\dvtx t\geq0\}$ in
$\mathsf{R}^{d}$ related to (\ref{logLap}). Recall also that for fixed
$ t>0$, with probability one, the measure state $ X_{t}%
$ is absolutely continuous (see \cite{Fleischmann1988critical}). The
following theorem is our \textit{main result}:
\begin{theorem}[(Dichotomy for densities)]\label{T.prop.dens}Fix $ t>0$
and $ X_{0}=\mu\in\mathcal{M}_{\mathrm{f} }$.

\begin{enumerate}[(a)]
\item[(a)] \textup{(Local H\"{o}lder continuity)}. If $ d=1$ and
$ \alpha>1+\beta$, then with probability one, there is a
continuous version $\tilde{X}_{t}$ of the density function of the measure
$ X_{t}(dx)$. Moreover, for each $ \eta<\eta
_{\mathrm{c}}:=\frac{\alpha}{1+\beta}-1$, this version $\tilde
{X}_{t}$ is locally H\"{o}lder continuous of index $ \eta$
\[
\sup_{x_{1},x_{2}\in K, x_{1}\neq x_{2}}\frac{|\tilde{X}_{t}%
(x_{1})-\tilde{X}_{t}(x_{2})|}{|x_{1}-x_{2}|^{\eta}} < \infty
\qquad\mbox{compact } K\subset\mathsf{R}.
\]

\item[(b)] \textup{(Optimal local H\"{o}lder index)}. Under conditions
as in the beginning of part \textup{(a)}, for every $ \eta\geq\eta_{\mathrm{c}}%
$ with probability one, for any open $ U\subseteq\mathsf{R}$,%
\[
\sup_{x_{1},x_{2}\in U, x_{1}\neq x_{2}}\frac{|\tilde{X}_{t}%
(x_{1})-\tilde{X}_{t}(x_{2})|}{|x_{1}-x_{2}|^{\eta}} = \infty
\qquad\mbox{whenever } X_{t}(U)>0.
\]

\item[(c)] \textup{(Local unboundedness)}. If $ d>1$ or $ \alpha
\leq1+\beta$, then with probability one, for all open
$ U\subseteq\mathsf{R}^{d}$,%
\[
\Vert X_{t}\Vert_{U} = \infty\qquad\mbox{whenever } X_{t}(U)>0.
\]
\end{enumerate}
\end{theorem}
\begin{remark}[(Any version)]
As in part (c), the statement in part (b) is valid also for any version
$X_{t}$ of the density function.
\end{remark}

\subsection{Some discussion}

At first sight, the result of Theorem \ref{T.prop.dens}(a), (b) is a bit
surprising. Let us recall again what is known about regularity
properties of
densities of $(\alpha,d,\beta)$-superprocesses. The case of continuous
super-Brownian motion ($\alpha=2$, $\beta=1$, $d=1)$ is very well
studied. As
already mentioned, densities exist at all times simultaneously, and
they are
locally H\"{o}lder continuous (in the spatial variable) for any index
$ \eta<\frac{1}{2}$. Moreover, it is known that
$\frac{1}{2}$ is optimal in this case. Now let us consider our result in
Theorem \ref{T.prop.dens}(a), (b), specialized to $\alpha=2$. Then we have
$ \eta_{\mathrm{c}}=\frac{2}{1+\beta}-1\downarrow0$ as
$ \beta\uparrow1$ where the limit 0 is different from the optimal
local H\"{o}lder index $\frac{1}{2}$ of continuous super-Brownian
motion. This
may confuse a reader and even raise a suspicion that something is wrong.
However there is an intuitive explanation for this discontinuity as we would
like to explain now.

Recall the notion of H\"{o}lder continuity \textit{at a point}. A
function $f$
is H\"{o}lder continuous with index $\eta\in(0,1)$ at a point $x_{0}$ if
there is a neighborhood $U(x_{0})$ such that
%
%
\begin{equation}
\vert f(x)-f(x_{0}) \vert\leq C |x-x_{0}|^{\eta
} \qquad\mbox{for all } x\in U(x_{0}).
\end{equation}
The \textit{optimal} H\"{o}lder index $H(x_{0})$ of $f$ at the point $x_{0}$
is defined as the supremum of all such $\eta$. Clearly, there are functions
where $H(x_{0})$ may vary with $x_{0 }$, and the index of a local
H\"{o}lder continuity in a domain cannot be larger than the smallest optimal
H\"{o}lder index at the points of the domain. The densities of continuous
super-Brownian motion are such that almost surely $H(x_{0})=\frac{1}{2}%
$ for all $x_{0 }$ whereas in our $\beta<1$ case of discontinuous
superprocesses the situation is quite different. The critical local H\"{o}lder
index $ \eta_{\mathrm{c}}=\frac{\alpha}{1+\beta}-1$ in our case is
a result of the influence of relatively high jumps of the superprocess that
occur close to time $ t$. So there are (random) points $x_{0}$
with $H(x_{0})=\eta_{\mathrm{c} }$. But these points are
\textit{exceptional} points; loosely speaking, there are not too many
of them.
We conjecture\footnote{We will verify this conjecture in an outcoming
extended version of \cite{FleischmannMytnikWachtel2009fixedWIASarxiv}.}
that at any \textit{given} point $x_{0}$ the optimal H\"{o}lder index
$H(x_{0}) $ equals $(\frac{1+\alpha}{1+\beta}-1)\wedge1=:\bar
{\eta}_{\mathrm{c}}>\eta_{\mathrm{c} }$. Now if $ \alpha
=2$, as $ \beta\uparrow1$ one gets the index
$ \frac{1}{2}$ corresponding to the case of continuous
super-Brownian motion.

This observation raises in fact a number of very interesting \textit{open
problems}:
\begin{conjecture}[(Multifractal spectrum)]
We conjecture that for any $ \eta\in(\eta_{\mathrm{c}},
\bar{\eta}_{\mathrm{c}})$ there are (random) points $x_{0}$ where the
density $X_{t}$ at the point $x_{0}$ is H\"{o}lder continuous with
index $\eta$. What is the \textit{Hausdorff dimension}, say $D(\eta)$, of
the (random) set $ \{x_{0}\dvtx H(x_{0})=\eta\}$? We
conjecture that%
%
%
\begin{equation}
\lim_{\eta\downarrow\eta_{\mathrm{c}}}D(\eta) = 0 \quad\mbox{and}\quad
\lim_{\eta\uparrow\bar{\eta}_{\mathrm{c}}}D(\eta)=1.
\end{equation}
This function $ \eta\mapsto D(\eta)$ reveals the so-called
\textit{multifractal} structure concerning the optimal H\"{o}lder index in
points for the densities of superprocesses with branching of index
$1+\beta<\alpha$ and is definitely worth studying. In
this connection, we refer to Jaffard \cite{Jaffard1999} where multifractal
properties of one-dimensional L\'{e}vy processes are studied.
\end{conjecture}

Another interesting direction would be a generalization of our results
to the
case of SPDEs driven by Levy noises. In recent years there has been increasing
interest in such SPDEs. Here we may mention the papers Saint Laubert Bi\'{e} \cite{Bie1998},
Mytnik \cite{Mytnik2002}, Mueller, Mytnik and Stan
\cite{MuellerMytnikStan2006} as well as Hausenblas \cite{Hausenblas2007}.
Note that in these papers properties of solutions are described in some
$\mathcal{L}^{p}$-sense. To the best of our knowledge not too many
things are
known about local H\"{o}lder continuity of solutions (in case of continuity).
The only result we know in this direction is \cite{MytnikPerkins2003} where
some local H\"{o}lder continuity of the fixed time density of super-Brownian
motion $(\alpha=2$, $\beta<1$, $d<\frac{2}{\beta}$, $a=0)$ was established.
However, the result there was far away from being optimal. With
Theorem \ref{T.prop.dens}(a), (b) we fill this gap. Our result also allows the
following conjecture:
\begin{conjecture}[(Regularity in case of SPDE with stable noise)]
\label{ConjEqu}%
Consider the SPDE,
%
%
\begin{equation}
\frac{\partial}{\partial t}X_{t}(x)=%
\bDelta_{\alpha}X_{t}(x)+g(X_{t-}(x))\dot{L}(t,x),
\end{equation}
where $\dot{L}$ is a $(1+\beta)$-stable noise without negative jumps,
and $g$
is such that solutions exist. Then there should exist versions of solutions
such that at fixed times regularity holds just as described in
Theorem \ref{T.prop.dens}(a), (b) with the same parameter classification,
in particular, with the same $\eta_{\mathrm{c} }$.
\end{conjecture}

\subsection{Martingale decomposition of $X$}

As in the $\alpha=2$ case of \cite{MytnikPerkins2003}, for the proof we need
the martingale decomposition of $ X$. For this purpose, we will
work with the following \textit{alternative description} of the
continuous-state branching mechanism $ \Psi$ from (\ref{not.Psi}):%
%
%
\begin{equation} \label{alt.Psi}%
\Psi(v) = -av + \varrho\int_{0}^{\infty} dr\, r^{-2-\beta} (
e^{-vr}-1+vr ) ,\qquad v\geq0,
\end{equation}
where $ $%
%
%
\begin{equation} \label{not.rho}%
\varrho:= b \frac{(1+\beta)\beta}{\Gamma(1-\beta)}
\end{equation}
with $ \Gamma$ denoting the famous Gamma function. The martingale
decomposition of $X$ in the following lemma is basically proven in
Dawson \cite{Dawson1993}, Section 6.1.

Denote by $\mathcal{C}_{\mathrm{b}}$ the set of all bounded and continuous
functions on $\mathsf{R}^{d}$. We add the sign $+$ if the functions are
additionally nonnegative. $\mathcal{C}_{\mathrm{b}}^{(k),+}$ with $k\geq1$
refers to the subset of functions which are $k$ times differentiable
and that
all derivatives up to the order $k$ belong to $\mathcal{C}_{\mathrm
{b}}^{+}$, too.
\begin{lemma}[(Martingale decomposition of $X$)]
\label{L.mart.dec}
Fix $ X_{0} = \mu\in\mathcal{M}_{\mathrm{f} }$.

\begin{enumerate}[(a)]
\item[(a)] \textup{(Discontinuities)}.
All discontinuities of the process $ X$
are jumps upward of the form $ r\delta_{x }$. More precisely, there
exists a random measure $ N (d(s,x,r) ) $ on $ \mathsf{R}_{+}\times
\mathsf{R}^{d}\times\mathsf{R}_{+}$ describing the jumps $ r\delta_{x}$
of $ X$ at times $ s$ at sites $x$ of size~$r$.

\item[(b)] \textup{(Jump intensities)}.
The compensator $ \hat{N}$ of
$ N$ is given by%
\[
\hat{N} (d(s,x,r) ) = \varrho
\,ds\, X_{s}(dx) r^{-2-\beta}\,dr;
\]
that is, $\tilde{N} := N-\hat{N}$ is a martingale
measure on $ \mathsf{R}_{+}\times\mathsf{R}^{d}\times\mathsf{R}_{+ }$.

\item[(c)] \textup{(Martingale decomposition)}.
For all $ \varphi\in\mathcal{C}_{\mathrm{b}}^{(2),+}$ and $ t\geq0$,%
\[
\langle X_{t},\varphi\rangle= \langle\mu,\varphi
\rangle+\int_{0}^{t} ds \langle X_{s},%
\bDelta_{\alpha}\varphi\rangle+M_{t}(\varphi)+a I_{t}(\varphi)
\]
with the discontinuous martingale%
\[
t \mapsto M_{t}(\varphi) := \int_{(0,t]\times\mathsf{R}^{d}\times
\mathsf{R}_{+}} \tilde{N} (d(s,x,r) )
r \varphi(x)
\]
and the increasing process%
\[
t \mapsto I_{t}(\varphi) := \int_{0}^{t} ds \langle
X_{s},\varphi\rangle.
\]
\end{enumerate}
\end{lemma}

From Lemma \ref{L.mart.dec} we get the related \textit{Green's function
representation},%
%
%
\begin{eqnarray}\label{Green}\quad
\langle X_{t},\varphi\rangle &=& \langle\mu
,S_{t}^{\alpha}\varphi\rangle+ \int_{(0,t]\times\mathsf{R}^{d}%
} M (d(s,x) ) S_{t-s}^{\alpha}%
\varphi(x)\nonumber\\[-8pt]\\[-8pt]
&&{} + a\int_{(0,t]\times\mathsf{R}^{d}} I
(d(s,x)) S_{t-s}^{\alpha}\varphi(x),\qquad
t\geq 0, \varphi\in\mathcal{C}_{\mathrm{b} }^{+},\nonumber
\end{eqnarray}
with $ M$ the martingale measure related to the martingale in part
(c) and $ I$ the measure related to the increasing process there.

We add also the following lemma which can be proved as Lemma 3.1 in Le Gall
and Mytnik \cite{LeGallMytnik2003}. For $p\geq1$, let $ \mathcal{L}%
_{\mathrm{loc}}^{p}(\mu)=\mathcal{L}_{\mathrm{loc}}^{p}(\mathsf{R}%
_{+}\times\mathsf{R}^{d}, S_{s}^{\alpha}\mu(x) \,ds
\,dx)$ denote the space of equivalence classes of measurable
functions $\psi$ such that%
%
%
\begin{equation}
\int_{0}^{T}ds\int_{\mathsf{R}^{d}}dx\, S_{s}^{\alpha}%
\mu(x) |\psi(s,x)|^{p} < \infty,\qquad T>0.
\end{equation}
\begin{lemma}[($L^{p}$-space with martingale measure)]\label{L.Lp}
Let $ X_{0}=\mu\in\mathcal{M}_{\mathrm{f} }$ and $ \psi\in\mathcal{L}
_{\mathrm{loc}}^{p}(\mu)$ for some
$ p\in(1+\beta,2)$. Then the martingale%
%
%
\begin{equation}
t\mapsto\int_{(0,t]\times\mathsf{R}^{d}} M (d
(s,x) ) \psi(s,x)
\end{equation}
is well defined.
\end{lemma}

Fix $ t>0$, $\mu\in\mathcal{M}_{\mathrm{f} }$. Suppose
$ d<\frac{\alpha}{\beta}$. Then the random measure $ X_{t}%
$ is a.s. absolutely continuous. From (\ref{Green}) we
get the following representation of a version of its \textit{density function}
(cf. \cite{MytnikPerkins2003,LeGallMytnik2003}):%
%
%
\begin{eqnarray}\label{rep.dens}%
X_{t}(x) & = & \mu\ast p_{t}^{\alpha} (x) + \int_{(0,t]\times
\mathsf{R}^{d}} M (d(s,y) )
p_{t-s}^{\alpha}(x-y)\nonumber\\
&&{} + a\int_{(0,t]\times\mathsf{R}^{d}} I
(d(s,y) ) p_{t-s}^{\alpha}(x-y)\\
&=&\!: Z_{t}^{1}(x)+Z_{t}^{2}(x)+Z_{t}^{3}(x),\qquad
x\in\mathsf{R}^{d},\nonumber
\end{eqnarray}
with notation in the obvious correspondence (and kernels $p^{\alpha}$
introduced in the beginning of Section \ref{SS.statement}).

This representation is the starting point for the proof of the local
H\"{o}lder continuity as claimed in Theorem \ref{T.prop.dens}(a). Main work
has to be done to deal with $Z_{t}^{2}$.

\subsection{Organization of the paper}

In Section \ref{sec:2} we develop some tools that will be used in the
following sections for the proof of Theorem \ref{T.prop.dens}. Also on
the way, in Section \ref{sec:2.3}, we are able to verify partially
Theorem \ref{T.prop.dens}(a) for some range of parameters $\alpha,\beta
$ using simple moment estimates. The proof of Theorem
\ref{T.prop.dens}(a) is completed in Section \ref{S.3} using a more
delicate analysis of the jump structure of the process. Section
\ref{sec:4} is devoted to the proof of part (c) of Theorem
\ref{T.prop.dens}. In Section \ref{sec:5}, which is the most
technically involved section, we verify Theorem \ref{T.prop.dens}(b).

\section{Auxiliary tools}
\label{sec:2}

In this section we always assume that $ d=1$.

\subsection{On the transition kernel of $\alpha$-stable motion}

The symbol $ C$ will always denote a generic positive constant,
which might change from place to place. On the other hand, $c_{(\#)}$ denotes
a constant appearing in formula line (or array) (\#).

We start with two estimates concerning the $\alpha$-stable transition kernel
$p^{\alpha}$.
\begin{lemma}[($\alpha$-stable density increment)]\label{L1}
For every $\delta\in[0,1]$,
%
%
\begin{equation} \label{L1.1}
|p_{t}^{\alpha}(x)-p_{t}^{\alpha}(y)| \leq C \frac{|x-y|^{\delta
}}{t^{\delta/\alpha}} \bigl(p_{t}^{\alpha}(x/2)+p_{t}^{\alpha}%
(y/2)\bigr),\qquad t>0,\ x,y\in\mathsf{R}.\hspace*{-29pt}
\end{equation}
\end{lemma}
\begin{pf}
For the case $\alpha=2$, see, for example, Rosen \cite{Rosen1987},
(2.4e). Suppose
$ \alpha<2$. It suffices to assume that $t=1$. In fact, multiply
$x,y$ by $t^{-1/\alpha}$ in the formula for the $t=1$ case, and use
that by
self-similarity, $p_{1}^{a}(t^{-1/a}x)=t^{1/\alpha}p_{t}^{\alpha}(x)$.

Now we use the well-known subordination formula
%
%
\begin{equation} \label{L1.2}%
p_{1}^{\alpha}(z) = \int_{0}^{\infty}ds\, q_{1}^{\alpha/2}%
(s) p_{s}^{(2)}(z),\qquad z\in\mathsf{R},
\end{equation}
where $q^{\alpha/2}$ denotes the continuous transition kernel of a stable
process on $\mathsf{R}_{+}$ of index $\alpha/2$, and by an abuse of notation,
$p^{(2)}$ refers to $p^{\alpha}$ in case $\alpha=2$. Consequently,
%
%
\begin{equation}
|p_{1}^{\alpha}(x)-p_{1}^{\alpha}(y)| \leq\int_{0}^{\infty
}ds\, q_{1}^{\alpha/2}(s) \bigl|p_{s}^{(2)}(x)-p_{s}^{(2)}(y)\bigr|.
\end{equation}
Hence, from the $ \alpha=2$ case,
%
%
\begin{eqnarray}
&&|p_{1}^{\alpha}(x)-p_{1}^{\alpha}(y)|\nonumber\\[-8pt]\\[-8pt]
&&\qquad \leq C |x-y|^{\delta}%
\int_{0}^{\infty}ds\, q_{1}^{\alpha/2}(s) s^{-\delta/2}%
\bigl(p_{s}^{(2)}(x/2)+p_{s}^{(2)}(y/2)\bigr).\nonumber
\end{eqnarray}
The lemma will be proved if we show that
%
%
\begin{equation} \label{L1.3}%
\int_{0}^{\infty}ds\, q_{1}^{\alpha/2}(s) s^{-\delta/2} p_{s}%
^{(2)}(x/2) \leq C p_{1}^{\alpha}(x/2),\qquad x\in\mathsf{R}.
\end{equation}
First, in view of (\ref{L1.2}),
%
%
\begin{equation} \label{L1.4}\qquad
\int_{1}^{\infty}ds\, q_{1}^{\alpha/2}(s) s^{-\delta/2} p_{s}%
^{(2)}(x/2) \leq\int_{1}^{\infty}ds\, q_{1}^{\alpha/2}%
(s) p_{s}^{(2)}(x/2) \leq p_{1}^{\alpha}(x/2).
\end{equation}
Second, by Brownian scaling,
%
%
\begin{eqnarray}\qquad
\int_{0}^{1}ds\, q_{1}^{\alpha/2}(s) s^{-\delta/2} p_{s}%
^{(2)}(x/2) &=& \int_{0}^{1}du\, q_{1}^{\alpha/2}(u) u^{-(\delta
+1)/2} p_{1}^{(2)} \biggl(\frac{x/2}{u^{1/2}} \biggr)\nonumber\\
&\leq& p_{1}^{(2)}(x/2)\int_{0}^{1}du\, q_{1}^{\alpha/2}%
(u) u^{-(\delta+1)/2}\\
&\leq& C p_{1}^{(2)}(x/2),\nonumber
\end{eqnarray}
where in the last step we have used the fact that $q_{1}^{\alpha/2}(u)$
decreases, as $u\downarrow0$, exponentially fast (cf.
\cite{Feller1971volII2nd}, Theorem 13.6.1). Since
$p_{1}^{(2)}(x/2)=\mathrm{o}(p_{1}%
^{\alpha}(x/2))$ as $x\uparrow\infty$, we have $p_{1}^{(2)}(x/2)\leq
Cp_{1}^{\alpha}(x/2), x\in\mathsf{R}$. Hence,
%
%
\begin{equation} \label{L1.5}%
\int_{0}^{1} ds\, q_{1}^{\alpha/2}(s) s^{-\delta/2} p_{s}%
^{(2)}(x/2)\leq C p_{1}^{\alpha}(x/2).
\end{equation}
Combining (\ref{L1.4}) and (\ref{L1.5}) gives (\ref{L1.3}), completing
the proof.
\end{pf}
\begin{lemma}[(Integrals of $\alpha$-stable density increment)]%
\label{L2}If $ \theta\in[1,1+\alpha)$ and $ \delta\in[
0,1]$ satisfy $ \delta<(1+\alpha-\theta)/\theta$ then%
%
%
\begin{eqnarray}\label{L2.1}\hspace*{28pt}
&& \int_{0}^{t}ds\int_{\mathsf{R}}dy\, p_{s}^{\alpha
}(y) |p_{t-s}^{\alpha}(x_{1}-y)-p_{t-s}^{\alpha}(x_{2}-y)|^{\theta}\nonumber\\[-8pt]\\[-8pt]
&&\qquad \leq C (1+t) |x_{1}-x_{2}|^{\delta\theta}\bigl(p_{t}^{\alpha}%
(x_{1}/2)+p_{t}^{\alpha}(x_{2}/2)\bigr),\qquad t>0, x_{1},x_{2}\in
\mathsf{R}.\nonumber
\end{eqnarray}
\end{lemma}
\begin{pf}
By Lemma \ref{L1}, for every $\delta\in[0,1]$,
%
%
\begin{eqnarray}
&& |p_{t-s}^{\alpha}(x_{1}-y)-p_{t-s}^{\alpha}(x_{2}-y)|^{\theta}
\nonumber\\[-8pt]\\[-8pt]
&&\qquad \leq C \frac{|x_{1}-x_{2}|^{\delta\theta}}{(t-s)^{\delta\theta/\alpha}}
\bigl(p_{t-s}^{\alpha}\bigl((x_{1}-y)/2\bigr)+p_{t-s}^{\alpha}%
\bigl((x_{2}-y)/2\bigr) \bigr)^{ \theta},\nonumber
\end{eqnarray}
$t>s\geq0$, $x_{1},x_{2},y\in\mathsf{R}$. Noting that
$p_{t-s}^{\alpha}(\cdot)\leq C (t-s)^{-1/\alpha}$, we obtain
%
%
\begin{eqnarray} \label{L2.2}\qquad
&&|p_{t-s}^{\alpha}(x_{1}-y)-p_{t-s}^{\alpha}(x_{2}-y)|^{\theta}\nonumber\\[-8pt]\\[-8pt]
&&\qquad \leq C \frac{|x_{1}-x_{2}|^{\delta\theta}}{(t-s)^{(\delta\theta
+\theta-1)/\alpha}} \bigl(p_{t-s}^{\alpha}\bigl((x_{1}-y)/2\bigr)+p_{t-s}%
^{\alpha}\bigl((x_{2}-y)/2\bigr) \bigr),\nonumber
\end{eqnarray}
$t>s\geq0$, $x_{1},x_{2},y\in\mathsf{R}$. Therefore,
\begin{eqnarray*}
&& \int_{0}^{t}ds\int_{\mathsf{R}}dy\, p_{s}^{\alpha
}(y) |p_{t-s}^{\alpha}(x_{1}-y)-p_{t-s}^{\alpha}(x_{2}-y)|^{\theta} \\
&&\qquad\leq C |x_{1}-x_{2}|^{\delta\theta}
\int_{0}^{t}ds (t-s)^{-(\delta\theta+\theta-1)/\alpha}\\
&&\qquad\quad{}\times\int_{\mathsf{R}}dy\, p_{s}^{\alpha}(y) \bigl(p_{t-s}^{\alpha
}\bigl((x_{1}-y)/2\bigr)+p_{t-s}^{\alpha}\bigl((x_{2}-y)/2\bigr) \bigr).
\end{eqnarray*}
By scaling of $ p^{\alpha}$,
%
%
\begin{eqnarray}
&&\int_{\mathsf{R}}dy\, p_{s}^{\alpha}(y) p_{t-s}^{\alpha
}\bigl((x-y)/2\bigr)\nonumber\\
&&\qquad= \frac{1}{2}\int_{\mathsf{R}}d%
y\, p_{2^{-\alpha}s}^{\alpha}(y/2) p_{t-s}^{\alpha}\bigl((x_{2}%
-y)/2\bigr) \nonumber\\
&&\qquad= \frac{1}{2} p_{2^{-\alpha}s+t-s}^{\alpha}(x/2) \\
&&\qquad= \frac{1}%
{2} (2^{-\alpha}s+t-s)^{-1/\alpha} p_{1}^{\alpha}\bigl((2^{-\alpha
}s+t-s)^{-1/\alpha}x/2\bigr)\nonumber\\
&&\qquad\leq t^{-1/\alpha} p_{1}^{\alpha}(t^{-1/\alpha}x/2) = p_{t}%
^{\alpha}(x/2),\nonumber
\end{eqnarray}
since $ 2^{-\alpha}t\leq2^{-\alpha}s+t-s\leq t$. As a result we
have the inequality
%
%
\begin{eqnarray}\qquad
&& \int_{0}^{t}ds\int_{\mathsf{R}}dy\, p_{s}^{\alpha
}(y) |p_{t-s}^{\alpha}(x_{1}-y)-p_{t-s}^{\alpha}(x_{2}-y)|^{\theta
}\nonumber\\[-8pt]\\[-8pt]
&&\qquad \leq C |x_{1}-x_{2}|^{\delta\theta}\bigl(p_{t}^{\alpha
}(x_{1}/2)+p_{t}^{\alpha}(x_{2}/2)\bigr)\int_{0}^{t}ds\, s^{-(\delta
\theta+\theta-1)/\alpha}.\nonumber
\end{eqnarray}
Noting that the latter integral is bounded by $ C (1+t)$, since
$ (\delta\theta+\theta-1)/\alpha<1$, we get the desired inequality.
\end{pf}

\subsection{An upper bound for a spectrally positive stable process}

Let $L=\{L_{t}\dvtx t\geq0\}$ denote a spectrally positive stable
process of
index $\kappa\in(1,2)$. Per definition, $L$ is an $\mathsf{R}$-valued
time-homogeneous process with independent increments and with Laplace
transform given by%
%
%
\begin{equation} \label{Laplace}%
\mathbf{E} e^{-\lambda L_{t}} = e^{t\lambda^{\kappa}%
},\qquad \lambda,t\geq0.
\end{equation}
Note that $L$ is the unique (in law) solution to the following martingale
problem:%
%
%
\begin{equation} \label{MP}%
t\mapsto e^{-\lambda L_{t}}-\int_{0}^{t} ds\, e
^{-\lambda L_{s}}\lambda^{\kappa} \mbox{ is a martingale for any }%
\lambda>0.
\end{equation}

Let $ \Delta L_{s}:=L_{s}-L_{s-}>0$ denote the jumps of $ L$.
\begin{lemma}[(Big values of the process in case of bounded
jumps)]\label{L3}
We have
%
%
\begin{eqnarray}
\mathbf{P} \Bigl( \sup_{0\leq u\leq t}L_{u}\mathsf{1}\Bigl\{\sup_{0\leq v\leq
u}\Delta L_{v}\leq y\Bigr\}\geq x \Bigr) \leq\biggl(\frac{C t}{xy^{\kappa
-1}} \biggr)^{ x/y},\nonumber\\[-8pt]\\[-8pt]
\eqntext{t>0, x,y>0.}
\end{eqnarray}
\end{lemma}
\begin{pf}
Since for $\tau>0$ fixed, $\{L_{\tau t}\dvtx t\geq0\}$ is equal to $\tau
^{1/\kappa}L$ in law, for the proof we may assume that $t=1$. Let $\{\xi
_{i}\dvtx i\geq1\}$ denote a family of independent copies of $L_{1 }$.
Set
%
%
\begin{equation}\quad
W_{ns} := \sum_{1\leq k\leq ns}\xi_{k},\qquad L_{s}^{(n)} := n^{-1/\kappa
}W_{ns},\qquad 0\leq s\leq1, n\geq1.
\end{equation}
Denote by $D_{[0,1]}$ the Skorohod space of c\`{a}dl\`{a}g functions
$ f\dvtx[0,1]\rightarrow\mathsf{R}$. For fixed $y>0$, let
$H\dvtx D_{[0,1]}\mapsto\mathsf{R}$ be defined by
%
%
\begin{equation}
H(f) = \sup_{0\leq u\leq1}f(u) \mathsf{1} \Bigl\{\sup_{0\leq v\leq u}\Delta
f(v)\leq y \Bigr\},\qquad f\in D_{[0,1] }.
\end{equation}
It is easy to verify that $H$ is continuous on the set
$D_{[0,1]}\setminus
J_{y}$ where $J_{y}:=\{f\in D_{[0,1]}\dvtx\Delta
f(v)=y\mbox{ for some }v\in[0,1]\}$. Since $\mathbf{P}(L\in
J_{y})=0$, from the invariance principle (see, e.g., Gikhman and Skorokhod
\cite{GikhmanSkorokhod1969}, Theorem 9.6.2) for $L^{(n)}$ we conclude
that%
%
%
\begin{equation}
\mathbf{P}\bigl(H(L)\geq x\bigr) = \lim_{n\uparrow\infty}\mathbf{P}%
\bigl(H\bigl(L^{(n)}\bigr)\geq x\bigr),\qquad x>0.
\end{equation}
Consequently, the lemma will be proved if we show that
%
%
\begin{eqnarray}\label{L3.1}%
\mathbf{P} \Bigl( \sup_{0\leq u\leq1}W_{nu}\mathsf{1}\Bigl\{\max_{1\leq k\leq
nu}\xi_{k}\leq yn^{1/\kappa}\Bigr\}\geq xn^{1/\kappa} \Bigr)
\leq\biggl(\frac{C}{xy^{\kappa-1}} \biggr)^{x/y},\nonumber\\[-8pt]\\[-8pt]
\eqntext{x,y>0, n\geq1.}
\end{eqnarray}
To this end, for fixed $ y^{\prime},h\geq0$, we consider the
sequence,%
%
%
\begin{equation}
\Lambda_{0}:=1,\qquad \Lambda_{n} := e^{hW_{n}}\mathsf{1}%
\Bigl\{\max_{1\leq k\leq n}\xi_{k}\leq y^{\prime}\Bigr\},\qquad n\geq1.
\end{equation}
It is easy to see that
%
%
\begin{equation}
\mathbf{E}\{\Lambda_{n+1} | \Lambda_{n}=e^{hu}\} = e
^{hu} \mathbf{E}\{e^{hL_{1}}; L_{1}\leq y^{\prime}\} \qquad
\mbox{for all } u\in\mathsf{R},
\end{equation}
and that
%
%
\begin{equation}
\mathbf{E}\{\Lambda_{n+1} | \Lambda_{n}=0\} = 0.
\end{equation}
In other words,
%
%
\begin{equation} \label{L3.2}%
\mathbf{E}\{\Lambda_{n+1} | \Lambda_{n}\} = \Lambda_{n} \mathbf{E}%
\{e^{hL_{1}}; L_{1}\leq y^{\prime}\}.
\end{equation}
This means that $\{\Lambda_{n}\dvtx n\geq1\}$ is a supermartingale
(submartingale) if $h$ satisfies $\mathbf{E}\{e^{hL_{1}}; L_{1}\leq
y^{\prime}\}\leq1 $ (respectively, $\mathbf{E}\{e^{hL_{1}} ;L_{1}\leq
y^{\prime}\}\geq1 $). If $\Lambda_{n}$ is a submartingale, then by
Doob's inequality,
%
%
\begin{equation}
\mathbf{P}\Bigl(\max_{1\leq k\leq n}\Lambda_{k}\geq e^{hx^{\prime}%
}\Bigr) \leq e^{-hx^{\prime}} \mathbf{E}\Lambda_{n},\qquad
x^{\prime}>0.
\end{equation}
But if $\Lambda_{n}$ is a supermartingale, then
%
%
\begin{equation}
\mathbf{P}\Bigl(\max_{1\leq k\leq n}\Lambda_{k}\geq e^{hx^{\prime}%
}\Bigr) \leq e^{-hx^{\prime}} \mathbf{E}\Lambda_{0 }%
= e^{-hx^{\prime}},\qquad x^{\prime}>0.
\end{equation}
From these inequalities and (\ref{L3.2}) we get
%
%
\begin{equation} \label{L3.3}%
\mathbf{P}\Bigl(\max_{1\leq k\leq n}\Lambda_{k}\geq e^{hx^{\prime}%
}\Bigr) \leq e^{-hx^{\prime}}\max\bigl\{1,(\mathbf{E}%
\{e^{hL_{1}}; L_{1}\leq y^{\prime}\})^{n} \bigr\}.
\end{equation}
It was proved by Fuk and Nagaev (\cite{FukNagaev71} see the first
formula in the proof of Theorem~4 there) that
\[
\mathbf{E}\{e^{hL_{1}}; L_{1}\leq y^{\prime}\} \leq1
+ h\mathbf{E}\{L_{1 }; L_{1}\leq y^{\prime}\}
+ \frac{e^{hy^{\prime}}-1-hy^{\prime}}{(y^{\prime})^{2}%
} V(y^{\prime}),\qquad h,y^{\prime}>0,
\]
where $ V(y^{\prime}):=\int_{-\infty}^{y^{\prime}}\mathbf{P}(L_{1}%
\in du) u^{2}>0$. Noting that the assumption
$\mathbf{E}L_{1}=0$ yields that $ \mathbf{E}\{L_{1 }; L_{1}\leq
y^{\prime
}\}\leq0$, we obtain
%
%
\begin{equation} \label{Nag}%
\mathbf{E}\{e^{hL_{1}}; L_{1}\leq y^{\prime}\} \leq1+\frac
{e^{hy^{\prime}}-1-hy^{\prime}}{(y^{\prime})^{2}} V(y^{\prime
}),\qquad h,y^{\prime}>0.
\end{equation}
Now note that%
%
%
\begin{eqnarray}\label{eq:2.30}
&&\Bigl\{\max_{1\leq k\leq n}W_{k}\mathsf{1}\Bigl\{\max_{1\leq i\leq k}\xi_{i}
\leq y^{\prime}\Bigr\}\geq x^{\prime} \Bigr\} \nonumber\\
&&\qquad= \Bigl\{\max_{1\leq k\leq
n}e^{hW_{k}} \mathsf{1}\Bigl\{\max_{1\leq i\leq k}\xi_{i}\leq y^{\prime
}\Bigr\}\geq e^{hx^{\prime}} \Bigr\}\\
&&\qquad = \Bigl\{\max_{1\leq k\leq n}\Lambda_{k}\geq e^{hx^{\prime}%
}\Bigr\}.\nonumber
\end{eqnarray}
%
Thus, combining (\ref{eq:2.30}), (\ref{Nag}) and (\ref{L3.3}), we get
%
%
\begin{eqnarray}
&&\mathbf{P} \Bigl(\max_{1\leq k\leq n}W_{k}\mathsf{1}\Bigl\{\max_{1\leq i\leq k}%
\xi_{i}\leq y^{\prime}\Bigr\}\geq x^{\prime} \Bigr)
\nonumber\\
&&\qquad\leq
\mathbf{P}\Bigl(\max_{1\leq k\leq n}\Lambda_{k}\geq e^{hx^{\prime}}\Bigr)\\
&&\qquad\leq\exp\biggl\{-hx^{\prime}+\frac{e^{hy^{\prime}}-1-hy^{\prime}%
}{(y^{\prime})^{2}} n V(y^{\prime}) \biggr\}.\nonumber
\end{eqnarray}
Choosing $h:=(y^{\prime})^{-1}\log(1+x^{\prime}y^{\prime}/n V(y^{\prime
}))$, we arrive, after some elementary calculations, at the bound,
%
%
\begin{equation}\hspace*{32pt}
\mathbf{P} \Bigl(\max_{1\leq k\leq n}W_{k}\mathsf{1}\Bigl\{\max_{1\leq i\leq k}%
\xi_{i}\leq y^{\prime}\Bigr\}\geq x^{\prime} \Bigr) \leq\biggl( \frac
{e n V(y^{\prime})}{x^{\prime}y^{\prime}} \biggr)^{ x^{\prime
}/y^{\prime}},\qquad x^{\prime},y^{\prime}>0.
\end{equation}
Since $\mathbf{P}(L_{1}>u)\sim C u^{-\kappa}$ as $u\uparrow\infty$, we have
$ V(y^{\prime})\leq C (y^{\prime})^{2-\kappa}$ for all $y^{\prime}>0$.
Therefore,
%
%
\begin{equation} \label{L3.4}\hspace*{32pt}
\mathbf{P} \Bigl(\max_{1\leq k\leq n}W_{k}\mathsf{1}\Bigl\{\max_{1\leq i\leq k}%
\xi_{i}\leq y^{\prime}\Bigr\}\geq x^{\prime} \Bigr) \leq\biggl(\frac{Cn}{x^{\prime
}(y^{\prime})^{\kappa-1}} \biggr)^{ x^{\prime}/y^{\prime}},\qquad x^{\prime
},y^{\prime}>0.
\end{equation}
Choosing finally $ x^{\prime}=xn^{1/\kappa}$, $y^{\prime
}=yn^{1/\kappa}$, we get (\ref{L3.1}) from (\ref{L3.4}). Thus, the
proof of the lemma is complete.
\end{pf}
\begin{lemma}[(Small process values)]\label{L.small.values}
There is a constant $ c_{\kappa}$ such that%
%
%
\begin{equation}
\mathbf{P} \Bigl(\inf_{u\leq t}L_{u}<-x \Bigr)\leq\exp\biggl\{-c_{\kappa}%
\frac{x^{\kappa/(\kappa-1)}}{t^{1/(\kappa-1)}} \biggr\},\qquad x,t>0.
\end{equation}
\end{lemma}
\begin{pf}
It is easy to see that for all $h>0$,%
%
%
\begin{equation}
\mathbf{P} \Bigl(\inf_{u\leq t}L_{u}<-x \Bigr)= \mathbf{P} \Bigl(\sup_{s\leq
t}e^{-hL_{u}}>e^{hx} \Bigr).
\end{equation}
Applying Doob's inequality to the submartingale $t\mapsto
e^{-hL_{t}}$, we obtain
%
%
\begin{equation}
\mathbf{P} \Bigl(\inf_{u\leq t}L_{u}<-x \Bigr)\leq e^{-hx}%
\mathbf{E} e^{-hL_{t}}.
\end{equation}
Taking into account definition (\ref{Laplace}), we have%
%
%
\begin{equation}
\mathbf{P} \Bigl(\inf_{u\leq t}L_{u}<-x \Bigr)\leq\exp\{-hx+th^{\kappa}\}.
\end{equation}
Minimizing the function $h\mapsto-hx+th^{\kappa}$, we get the
inequality in
the lemma with $c_{\kappa}=(\kappa-1)/(\kappa)^{\kappa/(\kappa-1)}$.
\end{pf}

\subsection{Local H\"{o}lder continuity with some index}\label{sec:2.3}

In this subsection we prove Theorem \ref{T.prop.dens}(a) for parameters
$\beta\geq\frac{\alpha-1}{2}$ (see Remark \ref{R.2.3}), whereas for parameters
$\beta<\frac{\alpha-1}{2}$ we obtain local H\"{o}lder continuity only with
nonoptimal bound on indexes. We use the Kolmogorov criterion for local
H\"{o}lder continuity to get these results. The proof of
Theorem \ref{T.prop.dens}(a) for parameters $\beta<\frac{\alpha-1}{2}$
will be finished in Section \ref{S.3}.

Fix $ t>0$, $\mu\in\mathcal{M}_{\mathrm{f} }$, and suppose $
\alpha>1+\beta$. Since our theorem is trivially valid for $\mu=0$, from
now on we everywhere suppose that $\mu\neq0$. Since we are
dealing with the case $ d=1$, the random measure $ X_{t}%
$ is a.s. absolutely continuous. Recall decomposition
(\ref{rep.dens}).

Clearly, the deterministic function $Z_{t}^{1}$ is Lipschitz
continuous by Lemma \ref{L1}. Next we turn to the random function $Z_{t}
^{3}$.
\begin{lemma}[(H\"{o}lder continuity of $Z_{t}^{3}$)]\label{L.HoeldZ3}
With probability one, $ Z_{t}^{3}$ is H\"{o}lder continuous of each
index $ \eta<\alpha-1$.
\end{lemma}
\begin{pf}
From Lemma \ref{L1} we get for fixed $\delta\in(0,\alpha-1)$,%
\[
|p_{t-s}^{\alpha}(x_{1}-y)-p_{t-s}^{\alpha}(x_{2}-y)| \leq
C \frac{|x_{1}-x_{2}|^{\delta}}{(t-s)^{(\delta+1)/\alpha}},
\qquad t>s>0, x_{1},x_{2},y\in\mathsf{R}.
\]
Therefore,
%
%
\begin{eqnarray}\label{2.38}%
&&|Z_{t}^{3}(x_{1})-Z_{t}^{3}(x_{2})|\nonumber\\
&&\qquad\leq |a|\int_{0}%
^{t} ds\int_{\mathsf{R}} X_{s}(dy) |p_{t-s}^{\alpha
}(x_{1}-y)-p_{t-s}^{\alpha}(x_{2}-y)|\nonumber\\[-8pt]\\[-8pt]
&&\qquad\leq C \Bigl(\sup_{s\leq t}X_{s}(\mathsf{R}) \Bigr) |x_{1}-x_{2}|^{\delta
}\int_{0}^{t} ds (t-s)^{-(\delta+1)/\alpha} \nonumber\\
&&\qquad\leq C \frac{\alpha}{\alpha-1-\delta} \Bigl(\sup_{s\leq t}X_{s}%
(\mathsf{R}) \Bigr) |x_{1}-x_{2}|^{\delta},\qquad
x_{1},x_{2}\in\mathsf{R}.\nonumber
\end{eqnarray}
Consequently,%
%
%
\begin{equation}
\sup_{x_{1}\neq x_{2}} \frac{|Z_{t}^{3}(x_{1})-Z_{t}^{3}(x_{2}%
)|}{|x_{1}-x_{2}|^{\delta}} < \infty\qquad\mbox{a.s.,}%
\end{equation}
and the proof is complete.
\end{pf}

Our main work concerns $Z_{t}^{2}$.
\begin{lemma}[($q$-norm)]\label{L4}
For each $\theta\in(1+\beta,2)$ and $q\in(1,1+\beta)$,%
%
%
\begin{eqnarray} \label{L4.1}%
&& \mathbf{E}|Z_{t}^{2}(x_{1})-Z_{t}^{2}(x_{2})|^{q}\nonumber\\
&&\qquad \leq C \biggl[ \biggl(\int_{0}^{t}ds\int_{\mathsf{R}}S_{s}^{\alpha
}\mu(dy) |p_{t-s}^{\alpha}(x_{1}-y)-p_{t-s}^{\alpha}%
(x_{2}-y)|^{\theta} \biggr)^{ q/\theta}\nonumber\\[-8pt]\\[-8pt]
&&\qquad\quad{} +\int_{0}^{t}ds\int_{\mathsf{R}}S_{s}^{\alpha}%
\mu(dy) |p_{t-s}^{\alpha}(x_{1}-y)-p_{t-s}^{\alpha}%
(x_{2}-y)|^{q} \biggr],\nonumber\\
\eqntext{x_{1},x_{2}\in\mathsf{R}.}
\end{eqnarray}
\end{lemma}

The proof can be done similarly to the proof of inequality (3.1) in
\cite{LeGallMytnik2003}.
\begin{corollary}[($q$-norm)]\label{L5}
For each $ \theta\in(1+\beta,2)$,
$q\in(1,1+\beta)$ and $ \delta>0$ satisfying $ \delta<\min
\{1,(1+\alpha-\theta)/\theta,(1+\alpha-q)/q\}$,
%
%
\begin{equation} \label{in.L5}%
\mathbf{E}|Z_{t}^{2}(x_{1})-Z_{t}^{2}(x_{2})|^{q} \leq
C |x_{1}-x_{2}|^{\delta q},\qquad x_{1},x_{2}\in\mathsf{R}.
\end{equation}
\end{corollary}
\begin{pf}
For every $\varepsilon\in(1,1+\alpha)$,
%
%
\begin{eqnarray}
&& \int_{0}^{t}ds\int_{\mathsf{R}}S_{s}^{\alpha}\mu(d%
y) |p_{t-s}^{\alpha}(x_{1}-y)-p_{t-s}^{\alpha}(x_{2}%
-y)|^{\varepsilon}\nonumber\\
&&\qquad = \int_{\mathsf{R}}\mu(dz)\int_{0}^{t}d%
s\int_{\mathsf{R}}dy\, p_{s}^{\alpha}(y-z) |p_{t-s}^{\alpha
}(x_{1}-z)-p_{t-s}^{\alpha}(x_{2}-z)|^{\varepsilon}\\
&&\qquad = \int_{\mathsf{R}}\mu(dz)\int_{0}^{t}d%
s\int_{\mathsf{R}}dy\, p_{s}^{\alpha}(y) |p_{t-s}^{\alpha}%
(x_{1}-z-y)-p_{t-s}^{\alpha}(x_{2}-z-y)|^{\varepsilon}.\hspace*{-25pt}\nonumber
\end{eqnarray}
Using Lemma \ref{L2}, we get for every positive $\delta<\min\{
1,(1+\alpha
-\varepsilon)/\varepsilon\}$,
\begin{eqnarray*}
&& \int_{0}^{t}ds\int_{\mathsf{R}}S_{s}^{\alpha}\mu(d%
y) |p_{t-s}^{\alpha}(x_{1}-y)-p_{t-s}^{\alpha}(x_{2}%
-y)|^{\varepsilon}\\
&&\qquad \leq C |x_{1}-x_{2}|^{\delta\varepsilon}\int_{\mathsf{R}}%
\mu(dz) \bigl(p_{t}^{\alpha}\bigl((x_{1}-z)/2\bigr)+p_{t}^{\alpha
}\bigl((x_{2}-z)/2\bigr) \bigr) \\
&&\qquad\leq C |x_{1}-x_{2}|^{\delta\varepsilon},
\end{eqnarray*}
since $ \mu,t$ are fixed. Applying this bound to both summands at
the right-hand side of (\ref{L4.1}) finishes the proof of the lemma.
\end{pf}
\begin{corollary}[(Finite $q$-norm of density)]\label{C1}
If $ K\subset\mathsf{R}$ is a compact and $ 1\leq q<1+\beta$, then
%
%
\begin{equation}
\mathbf{E} \Bigl(\sup_{x\in K}X_{t}(x) \Bigr)^{ q}<\infty.
\end{equation}
\end{corollary}
\begin{pf}
By Jensen's inequality, we may additionally assume that $q>1$. It
follows from
(\ref{rep.dens}) that%
%
%
\begin{equation}\hspace*{32pt}
\Bigl(\sup_{x\in K}X_{t}(x) \Bigr)^{ q} \leq4 \Bigl( \Bigl(\sup_{x\in K}%
\mu\ast p_{t}^{\alpha} (x) \Bigr)^{ q}+\sup_{x\in K}|Z_{t}%
^{2}(x)|^{q}+\sup_{x\in K}|Z_{t}^{3}(x)|^{q} \Bigr).
\end{equation}
Clearly, the first term at the right-hand side is finite. Furthermore,
according to Corollary 1.2 of Walsh \cite{Walsh1986}, inequality (\ref{in.L5})
implies that
%
%
\begin{equation}
\mathbf{E}\sup_{x\in K}|Z_{t}^{2}(x)|^{q}<\infty.
\end{equation}
Finally, proceeding as with the derivation of (\ref{2.38}), we obtain%
%
%
\begin{equation}
\sup_{x\in K}|Z_{t}^{3}(x)|\leq C \sup_{s\leq t}X_{s}%
(\mathsf{R}) \leq C e^{|a|t}\sup_{s\leq t} e^{-as}%
X_{s}(\mathsf{R}).
\end{equation}
Noting that $ s\mapsto e^{-as}X_{s}(\mathsf{R})$ is a
martingale, and using Doob's inequality, we conclude that%
%
%
\begin{equation}
\mathbf{E}\sup_{x\in K}|Z_{t}^{2}(x)|^{q} \leq C \mathbf{E}%
( e^{-at}X_{t}(\mathsf{R}))^{ q} < \infty.
\end{equation}
This completes the proof.
\end{pf}

Furthermore, Corollary \ref{L5} allows us to prove the following result:
\begin{proposition}[(Local H\"{o}lder continuity of $Z_{t}^{2}$)]\label{P1}
With probability one, $Z_{t}^{2}$ has a version which is locally
H\"{o}lder continuous of all orders $\eta>0$ satisfying
%
%
\begin{equation}\label{2cases}%
\eta< \eta_{\mathrm{c}}^{\prime} :=
\cases{\dfrac{\alpha}{1+\beta}-1, &\quad if
$\beta\geq(\alpha-1)/2$,\cr
\dfrac{\beta}{1+\beta}, &\quad if $\beta\leq(\alpha-1)/2$.}
\end{equation}
\end{proposition}
\begin{pf}
Let $\theta$, $q$ and $\delta$ satisfy the conditions in Corollary \ref{L5}.
Then almost surely $Z_{t}^{2}$ has a version which is locally H\"{o}lder
continuous of all orders smaller than $\delta-1/q$, (cf.
\cite{Walsh1986}, Corollary 1.2).

Let $\varepsilon>0$ satisfy $ \varepsilon<1-\beta$ and
$ \varepsilon<\beta$. Then $\theta=\theta_{\varepsilon}%
:=1+\beta+\varepsilon$ and $ q=q_{\varepsilon}:=1+\beta
-\varepsilon$ are in the range of parameters we are just
considering. Moreover, the condition $\delta<\min\{1,(1+\alpha
-\theta)/\theta,(1+\alpha-q)/q\}$ reads as%
%
%
\begin{equation}
\delta< \min\biggl\{1, \frac{\alpha-\beta-\varepsilon}{1+\beta+\varepsilon
} , \frac{\alpha-\beta+\varepsilon}{1+\beta-\varepsilon}%
\biggr\} =: f(\varepsilon).
\end{equation}
Hence, for all sufficiently small $\varepsilon>0$ we can choose $ \delta
=\delta_{\varepsilon}:=f(\varepsilon)-\varepsilon$. Thus,
$Z_{t}^{2}$ has a version which is locally H\"{o}lder continuous of all orders
smaller than $\delta_{\varepsilon}-1/q_{\varepsilon}$ for this
choice of $ \theta_{\varepsilon},q_{\varepsilon},\delta_{\varepsilon}%
$. Now%
\[
\delta_{\varepsilon}-\frac{1}{q_{\varepsilon}} \mathop{\longrightarrow
}_{\varepsilon
\downarrow0} \min\biggl\{1, \frac{\alpha-\beta}{1+\beta
} , \frac{\alpha-\beta}{1+\beta} \biggr\}-\frac{1}{1+\beta} = \min
\biggl\{1, \frac{\beta}{1+\beta} , \frac{\alpha-\beta-1}{1+\beta}
\biggr\},
\]
where this limit coincides with the claimed value of $ \eta_{\mathrm{c}%
}^{\prime}$, completing the proof.
\end{pf}
\begin{remark}[{[Proof of Theorem \ref{T.prop.dens}(a) for $\beta\geq
\frac{\alpha-1}{2}$]}]\label{R.2.3}
By Lemma \ref{L.HoeldZ3} and Proposition \ref{P1}, the
proof of Theorem \ref{T.prop.dens}(a) is finished for $\beta\geq\frac
{\alpha-1}{2}$.
\end{remark}

\subsection{Further estimates}

We continue to fix $t>0$, $\mu\in\mathcal{M}_{\mathrm{f}}%
\setminus\{0\}$, and to suppose $ \alpha>1+\beta$.
\begin{lemma}[(Local boundedness of uniformly smeared out
density)]\label{L6}
Fix a nonempty compact $ K\subset\mathsf{R}$ and a constant $c\geq1$.
Then%
%
%
\begin{equation}
V := V_{t}^{c}(K) := \sup_{0\leq s\leq t, x\in K}S_{c (t-s)}^{\alpha
}X_{s} (x) < \infty\qquad\mbox{almost surely}.
\end{equation}
\end{lemma}
\begin{pf}
Assume that the statement of the lemma does not hold,\break that is, there
exists an
event $A$ of positive probability such that\break $\sup_{0\leq s\leq t, x\in
K}S_{c (t-s)}^{\alpha}X_{s} (x)=\infty$ for every $\omega\in A$. Let
$n\geq1$. Put
\[
\tau_{n} := \cases{
\inf\bigl\{s<t\dvtx\mbox{there exists }x\in K\mbox{ such that }S_{c (t-s)}%
^{\alpha}X_{s} (x)>n \bigr\}, &\quad $\omega\in A$,\cr
t,&\quad $\omega\in A^{\mathrm{c}}$.}
\]
If $ \omega\in A$, choose $x_{n}=x_{n}(\omega)\in K$ such that
$ S_{c (t-\tau_{n})}^{\alpha}X_{\tau_{n}}(x_{n})>n$ whereas if
$ \omega\in A^{\mathrm{c}}$, take any $x_{n}=x_{n}(\omega)\in K$.
Using the strong Markov property gives
%
%
\begin{eqnarray}\label{56}
\mathbf{E}S_{(c-1)(t-\tau_{n})}^{\alpha}X_{t} (x_{n}) &=& \mathbf{EE}%
\bigl[S_{(c-1)(t-\tau_{n})}^{\alpha}X_{t} (x_{n}) | \mathcal{F}%
_{\tau_{n}}\bigr]\nonumber\\
&=& \mathbf{E} e^{a(t-\tau_{n})}S_{(c-1)(t-\tau_{n})}^{\alpha
}S_{(t-\tau_{n})}^{\alpha}X_{\tau_{n}}(x_{n}) \\
&\geq& e^{-|a|t}%
\mathbf{E}S_{c (t-\tau_{n})}^{\alpha}X_{\tau_{n}}(x_{n})\nonumber
\end{eqnarray}
[with $ e^{a(t-\tau_{n})}$ coming from the noncriticality of
branching in (\ref{not.Psi})]. From the definition of $(\tau
_{n},x_{n})$, we get%
%
%
\begin{equation} \label{to.infty}%
\mathbf{E}S_{c (t-\tau_{n})}^{\alpha}X_{\tau_{n}}(x_{n}) \geq n \mathbf
{P}%
(A)\rightarrow\infty\qquad\mbox{as } n\uparrow\infty.
\end{equation}
In order to get a contradiction, we want to prove boundedness in $n$ of the
expectation in (\ref{56}). If $c=1$, then%
%
%
\begin{equation}
\mathbf{E}X_{t}(x_{n}) \leq\mathbf{E}\sup_{x\in K}X_{t}(x) < \infty,
\end{equation}
the last step by Corollary \ref{C1}. Now suppose $c>1$. Choosing a compact
$K_{1}\supset K$ satisfying $\mathrm{dist}(K,(K_{1})^{\mathrm{c}%
})\geq1$, we have
\begin{eqnarray*}
&& \mathbf{E}S_{(c-1)(t-\tau_{n})}^{\alpha}X_{t} (x_{n})\\
&&\qquad = \mathbf{E}\int_{K_{1}}dy\, X_{t}(y) p_{(c-1)(t-\tau_{n}%
)}^{\alpha}(x_{n}-y)\\
&&\qquad\quad{} + \mathbf{E}\int_{(K_{1})^{\mathrm{c}}}d%
y\, X_{t}(y) p_{(c-1)(t-\tau_{n})}^{\alpha}(x_{n}-y)\\
&&\qquad \leq\mathbf{E}\sup_{y\in K_{1}}X_{t}(y)+\mathbf{E}X_{t}%
(\mathsf{R})\sup_{y\in(K_{1})^{\mathrm{c}}, x\in K, 0\leq s\leq
t }p_{(c-1)s}^{\alpha}(x-y).
\end{eqnarray*}
By our choice of $ K_{1}$ we obtain the bound,
%
%
\begin{equation}
\mathbf{E}S_{(c-1)(t-\tau_{n})}^{\alpha}X_{t} (x_{n}) \leq\mathbf{E}%
\sup_{y\in K_{1}}X_{t}(y)+C = C,
\end{equation}
the last step by Corollary \ref{C1}. Altogether, (\ref{56}) is bounded
in $n$,
and the proof is finished.
\end{pf}
\begin{lemma}[(Randomly weighted kernel increments)]\label{L7}
Fix $ \theta\in[1,1+\alpha)$, $ \delta\in[0,1]$ with
$ \delta<(1+\alpha-\theta)/\theta$, and a nonempty compact
$ K\subset\mathsf{R}$. Then
%
%
\begin{eqnarray}
&&\int_{0}^{t}ds\int_{\mathsf{R}}X_{s}(dy) |p_{t-s}%
^{\alpha}(x_{1}-y)-p_{t-s}^{\alpha}(x_{2}-y)|^{\theta}
\nonumber\\[-8pt]\\[-8pt]
&&\qquad\leq C V |x_{1}-x_{2}|^{\delta\theta},\qquad x_{1},x_{2}\in
K, \mbox{ a.s.},\nonumber
\end{eqnarray}
with $ V=V_{t}^{2^{\alpha}}(K)$ from Lemma \ref{L6}.
\end{lemma}
\begin{pf}
Using (\ref{L2.2}) gives
\begin{eqnarray*}
&& \int_{0}^{t}ds\int_{\mathsf{R}}X_{s}(dy) |p_{t-s}%
^{\alpha}(x_{1}-y)-p_{t-s}^{\alpha}(x_{2}-y)|^{\theta}\\
&&\qquad \leq
C |x_{1}-x_{2}|^{\delta\theta}
\int_{0}^{t}ds (t-s)^{-(\delta\theta+\theta-1)/\alpha}\\
&&\qquad\quad{}\times\int_{\mathsf{R}}X_{s}(dy) \bigl(p_{t-s}^{\alpha}\bigl((x_{1}%
-y)/2\bigr)+p_{t-s}^{\alpha}\bigl((x_{2}-y)/2\bigr) \bigr),
\end{eqnarray*}
uniformly in $ x_{1},x_{2}\in\mathsf{R}$. Recalling the scaling
property of $p^{\alpha}$, we get
\begin{eqnarray*}
&& \int_{0}^{t}ds\int_{\mathsf{R}}X_{s}(dy) |p_{t-s}%
^{\alpha}(x_{1}-y)-p_{t-s}^{\alpha}(x_{2}-y)|^{\theta}\\
&&\qquad \leq C |x_{1}-x_{2}|^{\delta\theta}\int_{0}^{t}d%
s (t-s)^{-(\delta\theta+\theta-1)/\alpha} \bigl(S_{2^{\alpha}(t-s)}^{\alpha
}X_{s}(x_{1})+S_{2^{\alpha}(t-s)}^{\alpha}X_{s}(x_{2}) \bigr).
\end{eqnarray*}
We complete the proof by applying Lemma \ref{L6}.
\end{pf}
\begin{remark}[(Lipschitz continuity of $Z_{t}^{3}$)]\label{R.Lipsch}
Using Lemma \ref{L7} with $\theta=1=\delta$, we see that $Z_{t}^{3}$ is
in fact a.s. Lipschitz continuous.
\end{remark}

Let $ \Delta X_{s}:=X_{s}-X_{s-}$ denote the jumps of the
measure-valued process $ X$.
\begin{lemma}[(Total jump mass)]\label{L8}
Let $\varepsilon>0$ and $ \gamma\in(0,(1+\beta)^{-1})$. There exists a
constant $c_{\mazinti{(\ref{inL8})}}=c_{\mazinti{(\ref{inL8})}}(\varepsilon,\gamma)$ such
that
%
%
\begin{equation} \label{inL8}%
\mathbf{P} \bigl(|\Delta X_{s}|>c_{\mazinti{(\ref{inL8})}} (t-s)^{(1+\beta)^{-1}-\gamma
}\mbox{ for some }s<t \bigr)\leq\varepsilon.
\end{equation}
\end{lemma}
\begin{pf}
Recall the random measure $N$ from Lemma \ref{L.mart.dec}(a). For any $c>0$,
set
%
%
\begin{eqnarray}\qquad
Y_{0} &:=& N \bigl([0,2^{-1}t)\times\mathsf{R}\times(c 2^{-\lambda
}t^{\lambda},\infty) \bigr),
\\
Y_{n} &:=& N
\bigl(\bigl[(1-2^{-n})t,(1-2^{-n-1})t\bigr)\nonumber\\[-8pt]\\[-8pt]
&&\hspace*{22.01pt}{}\times\mathsf{R}%
\times\bigl(c 2^{-\lambda(n+1)}t^{\lambda},\infty\bigr) \bigr),\qquad
n\geq1,\nonumber
\end{eqnarray}
where $\lambda:=(1+\beta)^{-1}-\gamma$. It is easy to see that
%
%
\begin{equation} \label{L12.1}\qquad
\mathbf{P} \bigl(|\Delta X_{s}|>c (t-s)^{\lambda}\mbox{ for some }%
s<t \bigr) \leq\mathbf{P} \Biggl(\sum_{n=0}^{\infty}Y_{n}\geq1 \Biggr) \leq
\sum_{n=0}^{\infty}\mathbf{E}Y_{n},
\end{equation}
where in the last step we have used the classical Markov inequality.
From the formula for the compensator $\hat{N}$ of $N$ in
Lemma \ref{L.mart.dec}(b),
%
%
\begin{equation}\quad
\mathbf{E}Y_{n} = \varrho\int_{(1-2^{-n})t}^{(1-2^{-n-1})t}d%
s\, \mathbf{E}X_{s}(\mathsf{R})\int_{c 2^{-\lambda(n+1)}t^{\lambda
}}^{\infty
}dr\, r^{-2-\beta},\qquad n\geq1.
\end{equation}
Now%
%
%
\begin{equation} \label{65}%
\mathbf{E}X_{s}(\mathsf{R}) = X_{0}(\mathsf{R}) e^{as} \leq
|\mu| e^{|a|t} =: c_{\mazinti{(\ref{65})}}.
\end{equation}
Consequently,
%
%
\begin{equation} \label{L12.2}%
\mathbf{E}Y_{n} \leq\frac{\varrho}{1+\beta} c_{\mazinti{(\ref{65})}}c^{-1-\beta
} 2^{-(n+1)\gamma(1+\beta)} t^{\gamma(1+\beta)}.
\end{equation}
Analogous calculations show that (\ref{L12.2}) remains valid also in
the case
$n=0$. Therefore,%
%
%
\begin{eqnarray}\label{L12.3}%
\sum_{n=0}^{\infty}\mathbf{E}Y_{n} &\leq& \frac{\varrho}{1+\beta}
c_{\mazinti{(\ref{65})}}c^{-1-\beta} t^{\gamma(1+\beta)}\sum_{n=0}^{\infty}
2^{-(n+1)\gamma(1+\beta)}\nonumber\\[-8pt]\\[-8pt]
&=& \frac{\varrho}{1+\beta} c_{\mazinti{(\ref{65})}}c^{-1-\beta} t^{\gamma
(1+\beta)}
\frac{2^{-\gamma(1+\beta)}}{1-2^{-\gamma(1+\beta)}}.\nonumber
\end{eqnarray}
Choosing $c=c_{\mazinti{(\ref{inL8})}}$ such that the expression in (\ref{L12.3}) equals
$\varepsilon$, and combining with (\ref{L12.1}), the proof is complete.
\end{pf}

\subsection{Representation as time-changed stable process}

We return to general \mbox{$t>0$}. Recall the martingale measure $ M$
related to the martingale in Lemmas \ref{L.mart.dec}(c) and \ref{L.Lp}.
\begin{lemma}[(Representation as time-changed stable process)]\label{L9}
Suppose $ p\in(1+\beta,2)$ and let $ \psi\in\mathcal{L}_{\mathrm
{loc}}^{p}(\mu)$
with $\psi\geq0$. Then there exists a spectrally positive $(1+\beta)$-stable
process $\{L_{t}\dvtx t\geq0\}$ such that%
%
%
\begin{equation}
Z_{t}(\psi) := \int_{(0,t]\times\mathsf{R}} M
(d(s,y) ) \psi(s,y) = L_{T(t) },\qquad t\geq0,
\end{equation}
where $T(t):=\int_{0}^{t}d%
s\int_{\mathsf{R}}X_{s}(dy) (\psi(s,y))^{1+\beta}$.
\end{lemma}
\begin{pf}
Let us write It\^{o}'s formula for $ e^{-Z_{t}(\psi)}$%
%
%
\begin{eqnarray}
e^{-Z_{t}(\psi)}-1
&=& \mbox{local martingale}\nonumber\\
&&{} +\varrho\int_{0}^{t}ds\, e^{-Z_{s}(\psi)}%
\int_{\mathsf{R}}X_{s}(dy)\\
&&\hspace*{10pt}{}\times\int_{0}^{\infty}d%
r \bigl( e^{-r\psi(s,y)}-1+r \psi(s,y)\bigr) r^{-2-\beta}.\nonumber
\end{eqnarray}
Define $ \tau(t):=T^{-1}(t)$, and put $ t^{\ast}:=\inf\{
t\dvtx\tau(t)=\infty\}$. Then it is easy to get for every
$v>0$,
%
%
\begin{eqnarray}\hspace*{32pt}
e^{-vZ_{\tau(t)}(\psi)} & = & 1+\int_{0}^{t}d%
s\, e^{-vZ_{\tau(s)}(\psi)} \frac{X_{\tau(s)}(v^{1+\beta}%
\psi^{1+\beta}(s,\cdot))}{X_{\tau(s)}(\psi^{1+\beta}(s,\cdot
))}+\mbox{loc. mart.}\nonumber\\[-8pt]\\[-8pt]
& = & 1+\int_{0}^{t}ds\, e^{-vZ_{\tau(s)}(\psi)} v^{1+\beta
}+\mbox{loc. mart.,}\qquad t\leq t^{\ast}.\nonumber
\end{eqnarray}
Since the local martingale is bounded, it is in fact a martingale. Let
$ \tilde{L}$ denote a spectrally positive process of index
$1+\beta$, independent of $X$. Define%
%
%
\begin{equation}
L_{t} := \cases{Z_{\tau(t)}(\psi), &\quad $t\leq t^{\ast}$,\cr
Z_{\tau(t^{\ast})}(\psi)+\tilde{L}_{t-t^{\ast} }, &\quad $t>t^{\ast}$
\mbox{ (if $t^{\ast}<\infty$)}.}
\end{equation}
Then we can easily get that $L$ satisfies the martingale problem (\ref{MP})
with $\kappa$ replaced by $1+\beta$. Now by time change back we obtain%
%
%
\begin{equation}
Z_{t}(\psi)=\tilde{L}_{T(t)}=L_{T(t) },
\end{equation}
completing the proof.
\end{pf}

\section{Local H\"{o}lder continuity}\label{S.3}

\mbox{}

\begin{pf*}{Proof of Theorem \protect\ref{T.prop.dens}\textup{(a)}}
We continue to assume that $d=1$, and that $ t>0$ and
$\mu\in\mathcal{M}_{\mathrm{f}}\setminus\{0\}$ are fixed. For $\beta
\geq(\alpha-1)/2$ the desired existence of a locally H\"{o}lder continuous
version of $ Z_{t}^{2}$ of required orders is already proved in
Proposition \ref{P1}. Therefore, in what follows we shall consider the
complementary case $\beta<(\alpha-1)/2$. Fix any compact set $K$ and
$x_{1}<x_{2}$ belonging to it. By definition (\ref{rep.dens}) of
$Z_{t}^{2}%
$,
%
%
\begin{eqnarray}\quad
Z_{t}^{2}(x_{1})-Z_{t}^{2}(x_{2}) &=& \int_{(0,t]\times\mathsf{R}} M (
d(s,y) ) \bigl(p_{t-s}^{\alpha}(x_{1}%
-y)-p_{t-s}^{\alpha}(x_{2}-y)\bigr)\nonumber\\
&=& \int_{(0,t]\times\mathsf{R}} M (d%
(s,y) ) \varphi_{+}(s,y)\\
&&{}-\int_{(0,t]\times\mathsf{R}} M
(d(s,y)) \varphi_{-}(s,y),\nonumber
\end{eqnarray}
where $\varphi_{+}(s,y)$ and $\varphi_{-}(s,y)$ are the positive and negative
parts of $ p_{t-s}^{\alpha}(x_{1}-y)-p_{t-s}^{\alpha}(x_{2}-y)$. It is easy
to check that $\varphi_{+}$ and $\varphi_{-}$ satisfy the assumptions in
Lem\-ma~\ref{L9}. Thus, there exist stable processes $L^{1}$ and $L^{2}$ such
that
%
%
\begin{equation}\label{T1.1}%
Z_{t}^{2}(x_{1})-Z_{t}^{2}(x_{2}) = L_{T_{+}}^{1}-L_{T_{- }}^{2},
\end{equation}
where $ T_{\pm}:=\int_{0}^{t}d%
s\int_{\mathsf{R}}X_{s}(dy) (\varphi_{\pm}(s,y))^{1+\beta
}$.

The idea behind the proof of the existence of the required version of
$ Z_{t}^{2}$ is as follows. We first control the jumps of $L^{1}$ and $L^{2}$
for $t\leq T_{\pm}$ and then use Lem\-ma~\ref{L3} to get the necessary bounds
on $L_{T_{+ }}^{1},L_{T_{- }}^{2}$ themselves.

Fix any $\varepsilon\in(0,1)$. According to Lemma \ref{L6}, there
exists a constant $c_{\varepsilon}$ such that
%
%
\begin{equation}\label{13_04}%
\mathbf{P}(V\leq c_{\varepsilon})\geq1-\varepsilon,
\end{equation}
where $ V=V_{t}^{2^{\alpha}}(K)$. Consider again $ \gamma
\in(0,(1+\beta)^{-1})$ and set
%
%
\begin{equation} \label{A.eps}%
A^{\varepsilon} := \bigl\{|\Delta X_{s}|\leq c_{\mazinti{(\ref{inL8})}}%
(t-s)^{(1+\beta)^{-1}-\gamma}\mbox{ for all }s<t \bigr\}\cap\{V\leq
c_{\varepsilon}\}.
\end{equation}
By Lemma \ref{L8} and by (\ref{13_04}),
%
%
\begin{equation} \label{A.eps.est}%
\mathbf{P}(A^{\varepsilon})\geq1-2\varepsilon.
\end{equation}

Define $Z_{t}^{2,\varepsilon}(x):=Z_{t}^{2}(x)\mathsf{1}(A^{\varepsilon})$.
We first show that $Z_{t}^{2,\varepsilon}$ has a version which
is locally H\"{o}lder continuous of all orders $\eta$ smaller than
$\eta_{\mathrm{c} }$. It follows from (\ref{T1.1}) that
%
%
\begin{eqnarray} \label{T1.2}%
&& \mathbf{P} \bigl(|Z_{t}^{2,\varepsilon}(x_{1})-Z_{t}^{2,\varepsilon
}(x_{2})|\geq2r |x_{1}-x_{2}|^{\eta} \bigr)\nonumber\\
&&\qquad \leq\mathbf{P}(L_{T_{+}}^{1}\geq r |x_{1}-x_{2}|^{\eta
}, A^{\varepsilon})\\
&&\qquad\quad{}+\mathbf{P}(L_{T_{-}}^{2}\geq r |x_{1}%
-x_{2}|^{\eta}, A^{\varepsilon}),\qquad r>0.\nonumber
\end{eqnarray}
Note that on $A^{\varepsilon}$ the jumps of $M (
d(s,y) ) $ do not exceed
$c_{\mazinti{(\ref{inL8})} }(t-s)^{(1+\beta)^{-1}-\gamma}$ since the jumps
of $X$ are bounded by the same values on $A^{\varepsilon}$. Hence the
jumps of
the process $ u\mapsto$ $\int_{(0,u]\times\mathsf{R}} M (
d(s,y) ) \varphi_{\pm}(s,y)$ are bounded by%
%
%
\begin{equation} \label{T1.3}%
c_{\mazinti{(\ref{inL8})}} \sup_{s<t}(t-s)^{(1+\beta)^{-1}-\gamma} \sup_{y\in
\mathsf{R}}\varphi_{\pm}(s,y).
\end{equation}
Obviously,
%
%
\begin{equation} \label{T1.4}%
\sup_{y\in\mathsf{R}}\varphi_{\pm}(s,y) \leq{\sup_{y\in\mathsf{R}%
}} |p_{t-s}^{\alpha}(x_{1}-y)-p_{t-s}^{\alpha}(x_{2}-y) |.
\end{equation}
Assume additionally that $ \gamma<\eta_{\mathrm{c}}/\alpha$. Using
Lemma \ref{L1} with $\delta=\eta_{\mathrm{c}}-\alpha\gamma$ gives
%
%
\begin{eqnarray} \label{T1.5}%
&& {\sup_{y\in\mathsf{R}}} |p_{t-s}^{\alpha}(x_{1}-y)-p_{t-s}^{\alpha}%
(x_{2}-y) | \nonumber\\
&&\qquad \leq C |x_{1}-x_{2}|^{\eta_{\mathrm{c}}-\alpha\gamma} (t-s)^{-\eta
_{\mathrm{c}}/\alpha+\gamma} \sup_{z\in\mathsf{R}}p_{t-s}^{\alpha
}(z)\nonumber\\[-8pt]\\[-8pt]
&&\qquad \leq C |x_{1}-x_{2}|^{\eta_{\mathrm{c}}-\alpha\gamma} (t-s)^{-\eta
_{\mathrm{c}}/\alpha+\gamma} (t-s)^{-1/\alpha}\nonumber\\
&&\qquad = C |x_{1}-x_{2}|^{\eta_{\mathrm{c}}-\alpha\gamma} (t-s)^{-
{1}/({1+\beta})+\gamma}.\nonumber
\end{eqnarray}
Combining (\ref{T1.3})--(\ref{T1.5}), we see that all jumps of
$ u\mapsto\int_{(0,u]\times\mathsf{R}} M (
d(s,y) ) \varphi_{\pm}(s,y)$ on the set $A^{\varepsilon}$ are
bounded by
%
%
\begin{equation} \label{C}%
c_{\mazinti{(\ref{C})}} |x_{1}-x_{2}|^{\eta_{\mathrm{c}}-\alpha\gamma}
\end{equation}
for some constant $ c_{\mazinti{(\ref{C})}}=c_{\mazinti{(\ref{C})}}(\varepsilon)$.
Therefore, by an abuse of notation writing $L_{T_{\pm}}$ for $L_{T_{+}}^{1}$
and $L_{T_{- }}^{2}$,
%
%
\begin{eqnarray}\hspace*{8pt}
&& \mathbf{P}(L_{T_{\pm}}\geq r |x_{1}-x_{2}|^{\eta}, A^{\varepsilon
}) \nonumber\\
&&\qquad = \mathbf{P} \Bigl(L_{T_{\pm}}\geq r |x_{1}-x_{2}|^{\eta}, \sup
_{u<T_{\pm}}\Delta L_{u}\leq c_{\mazinti{(\ref{C})}} |x_{1}-x_{2}|^{\eta_{\mathrm
{c}%
}-\alpha\gamma}, A^{\varepsilon} \Bigr)\\
&&\qquad \leq\mathbf{P} \Bigl( \sup_{v\leq T_{\pm}}L_{v}\mathsf{1} \Bigl\{\sup
_{u<v}\Delta L_{u}\leq c_{\mazinti{(\ref{C})}} |x_{1}-x_{2}|^{\eta_{\mathrm{c}}%
-\alpha\gamma} \Bigr\}\geq r |x_{1}-x_{2}|^{\eta}, A^{\varepsilon} \Bigr).\nonumber\hspace*{-12pt}
\end{eqnarray}
Since
%
%
\begin{equation}
T_{\pm} \leq\int_{0}^{t}ds\int_{\mathsf{R}}X_{s}(d%
y) |p_{t-s}^{\alpha}(x_{1}-y)-p_{t-s}^{\alpha}(x_{2}-y) |^{1+\beta},
\end{equation}
applying Lemma \ref{L7} with $\theta=1+\beta$ and $\delta=1$ (since
$\beta<(\alpha-1)/2$), we get the bound
%
%
\begin{equation} \label{81}%
T_{\pm} \leq c_{\mazinti{(\ref{81})}} |x_{1}-x_{2}|^{1+\beta} \qquad\mbox{on } \{V\leq
c_{\varepsilon}\},
\end{equation}
for some $c_{\mazinti{(\ref{81})}}=c_{\mazinti{(\ref{81})}}(\varepsilon)$. Consequently,
\begin{eqnarray*}
&& \mathbf{P}(L_{T_{\pm}}\geq r |x_{1}-x_{2}|^{\eta}, A^{\varepsilon
})\\
&&\qquad \leq\mathbf{P} \Bigl( \sup_{v\leq c_{\mazinti{(\ref{81})}}|x_{1}-x_{2}|^{1+\beta}%
}L_{v} \mathsf{1}\Bigl\{\sup_{u<v}\Delta L_{u}\leq c_{\mazinti{(\ref{C})}}%
|x_{1}-x_{2}|^{\eta_{\mathrm{c}}-\alpha\gamma}\Bigr\}\\
&&\hspace*{232.3pt}\geq r |x_{1}%
-x_{2}|^{\eta} \Bigr).
\end{eqnarray*}
Using Lemma \ref{L3} with $ \kappa=1+\beta$, $t=c_{\mazinti{(\ref{81}%
)}}|x_{1}-x_{2}|^{1+\beta}$, $x=r |x_{1}-x_{2}|^{\eta}$,
and $ y=c_{\mazinti{(\ref{C})}}|x_{1}-x_{2}|^{\eta_{\mathrm{c}%
}-\alpha\gamma}$, and noting that%
%
%
\begin{eqnarray}
1+\beta-\eta-\beta(\eta_{\mathrm{c}}-\alpha\gamma) &=& 2+2\beta-\alpha
+(\eta_{\mathrm{c}}-\eta)+\beta\alpha\gamma\nonumber\\[-8pt]\\[-8pt]
&>& 2+2\beta-\alpha,\nonumber
\end{eqnarray}
we obtain
%
%
\begin{eqnarray} \label{3.14}%
&& \mathbf{P}(L_{T_{\pm}}\geq r |x_{1}-x_{2}|^{\eta}, A^{\varepsilon
}) \nonumber\\[-8pt]\\[-8pt]
&&\qquad \leq\bigl(c_{\mazinti{(\ref{3.14})}} r^{-1}|x_{1}-x_{2}|^{(2\beta+2-\alpha
)}\bigr)^{(c_{\mmazinti{(\ref{C})}}^{-1}r|x_{1}-x_{2}|^{\eta-\eta_{\mathrm{c}%
}+\alpha\gamma})}\nonumber
\end{eqnarray}
for some $c_{\mazinti{(\ref{3.14})}}=c_{\mazinti{(\ref{3.14})}}(\varepsilon)$. Applying
this bound with $\gamma=(\eta_{\mathrm{c}}-\eta)/2\alpha$ to the
summands at the right-hand side in (\ref{T1.2}), and noting that $
2\beta+2-\alpha$ is also constant here, we have
%
%
\begin{eqnarray} \label{3.16}%
&& \mathbf{P} \bigl(|Z_{t}^{2,\varepsilon}(x_{1})-Z_{t}^{2,\varepsilon
}(x_{2})|\geq2r |x_{1}-x_{2}|^{\eta} \bigr) \nonumber\\[-8pt]\\[-8pt]
&&\qquad \leq2\bigl(c_{\mazinti{(\ref{3.14})}} r^{-1}|x_{1}-x_{2}%
|\bigr)^{(c_{\mmazinti{(\ref{3.16})}}r|x_{1}-x_{2}|^{(\eta-\eta_{\mathrm{c}}%
)/2})}.\nonumber
\end{eqnarray}
This inequality yields that all the conditions of Theorem III.5.6 of\break Gihman
and Skorokhod \cite{GikhmanSkorokhod1974} hold with $g(h)=2h^{\eta}$ and
$q(r,h)=\break 2(c_{\mazinti{(\ref{3.14})}} r^{-1}h)^{(c_{\mmazinti{(\ref{3.16}%
)}}rh^{(\eta-\eta_{\mathrm{c}})/2})}$, from which we conclude that almost
surely $Z_{t}^{2,\varepsilon}$ has a version which is locally H\"{o}lder
continuous of all orders $\eta<\eta_{\mathrm{c} }$.

By an abuse of notation, from now on the symbol $Z_{t}^{2,\varepsilon}$ always
refers to this continuous version. Consequently,
%
%
\begin{equation} \label{T1.6}%
\lim_{k\uparrow\infty}\mathbf{P} \biggl( \sup_{x_{1},x_{2}\in K, x_{1}\neq
x_{2}}\frac{|Z_{t}^{2,\varepsilon}(x_{1})-Z_{t}^{2,\varepsilon}%
(x_{2})|}{|x_{1}-x_{2}|^{\eta}}>k \biggr) = 0.
\end{equation}
Combining this with the bound
%
%
\begin{eqnarray}\quad
&& \mathbf{P} \biggl( \sup_{x_{1},x_{2}\in K, x_{1}\neq x_{2}}\frac
{|Z_{t}^{2}(x_{1})-Z_{t}^{2}(x_{2})|}{|x_{1}-x_{2}|^{\eta}%
}>k \biggr) \nonumber\\[-8pt]\\[-8pt]
&&\qquad \leq\mathbf{P} \biggl( \sup_{x_{1},x_{2}\in K, x_{1}\neq x_{2}}%
\frac{|Z_{t}^{2,\varepsilon}(x_{1})-Z_{t}^{2,\varepsilon}(x_{2}%
)|}{|x_{1}-x_{2}|^{\eta}}>k, A^{\varepsilon} \biggr)+\mathbf{P}%
(A^{\varepsilon,\mathrm{c}})\nonumber
\end{eqnarray}
(with $A^{\varepsilon,\mathrm{c}}$ denoting the complement of
$A^{\varepsilon
})$, gives
%
%
\begin{equation}
\limsup_{k\uparrow\infty}\mathbf{P} \biggl( \sup_{x_{1},x_{2}\in K, x_{1}\neq
x_{2}}\frac{|Z_{t}^{2}(x_{1})-Z_{t}^{2}(x_{2})|}{|x_{1}-x_{2}%
|^{\eta}}>k \biggr) \leq2\varepsilon.
\end{equation}
Since $\varepsilon$ may be arbitrarily small, this immediately implies
%
%
\begin{equation}
\sup_{x_{1},x_{2}\in K, x_{1}\neq x_{2}}\frac{|Z_{t}^{2}(x_{1}%
)-Z_{t}^{2}(x_{2})|}{|x_{1}-x_{2}|^{\eta}}<\infty, \qquad\mbox{almost
surely}.
\end{equation}
This is the desired local H\"{o}lder continuity of $ Z_{t}^{2}$,
for all $\eta<\eta_{\mathrm{c} }$. Because
$ \eta_{\mathrm{c}}<\alpha-1$, together with Lemma \ref{L.HoeldZ3}
the proof of Theorem \ref{T.prop.dens}(a) is complete.
\end{pf*}

\section[Local unboundedness: Proof of Theorem 1.2(c)]{Local unboundedness: Proof of
Theorem \protect\ref{T.prop.dens}\textup{(c)}}
\label{sec:4}

In the proof we use ideas from the proofs of Theorems 1.1(b) and 1.2 of
\cite{MytnikPerkins2003}. Throughout this section, suppose $d>1$ or
$\alpha\leq1+\beta$. Recall that $ t>0$ and
$ X_{0}=\mu\in\mathcal{M}_{\mathrm{f}}\setminus\{0\}$ are fixed.
We want to verify that for each version of the density function
$X_{t}$ the
property
%
%
\begin{equation} \label{basicprop}
\Vert X_{t} \Vert_{B}=\infty\qquad\mathbf{P}\mbox{-a.s. on the event
} \{X_{t}(B)>0\}%
\end{equation}
holds whenever $B$ is a fixed open ball in $ \mathsf{R}^{d}$. Then
the claim of Theorem \ref{T.prop.dens}(c) follows as in the proof of
Theorem 1.1(b) in \cite{MytnikPerkins2003}. We thus fix such $B$.

As in \cite{MytnikPerkins2003} to get (\ref{basicprop}) we first show
that on the event $\{X_{t}(B)>0\}$ there are always sufficiently
``big'' jumps of $X$ that occur close to time $t$. This is done in
Lemma \ref{lem:3} below. Then with the help of properties of the
log-Laplace equation derived in Lemma \ref{lem:2} we are able to show
that the ``big'' jumps are large enough to ensure the unboundedness of
the density at time $t$. Loosely speaking the density is getting
unbounded in the proximity of big jumps.

In order to fulfil the above program, we start with deriving the
continuity of
$X_{\cdot}(B)$ at (fixed) time $t$.
\begin{lemma}[(Path continuity at fixed times)]\label{lem:1}
For the fixed $t>0$,%
%
%
\begin{equation}
\lim_{s\rightarrow t}X_{s}(B) = X_{t}(B) \qquad\mbox{a.s.}%
\end{equation}
\end{lemma}
\begin{pf}
Since $t$ is fixed, $X$ is continuous at $t$ with probability $1$. Therefore,
%
%
\begin{equation}
X_{t}(B) \leq\liminf_{s\rightarrow t}X_{s}(B) \leq\limsup
_{s\rightarrow
t}X_{s}(B) \leq\limsup_{s\rightarrow t}X_{s}(\overline{B}) \leq
X_{t}(\overline{B})
\end{equation}
with $\overline{B}$ denoting the closure of $B$. But since $X_{t}%
(dx)$ is absolutely continuous with respect to Lebesgue measure, we
have $X_{t}(B)=X_{t}(\overline{B})$. Thus the proof is complete.
\end{pf}
\begin{lemma}[(Explosion)]\label{L.explosion}
Let $f\dvtx(0,t)\rightarrow(0,\infty)$ be
measurable such that
%
%
\begin{equation}
\int_{t-\delta}^{t}ds\, f(t-s)=\infty\qquad\mbox{for all sufficiently
small } \delta\in(0,t).
\end{equation}
Then for these $\delta$,
%
%
\begin{equation}
\quad \int_{t-\delta}^{t}ds\, X_{s}(B)f(t-s) = \infty,\qquad\mathbf{P}%
\mbox{-a.s. on the event }\{X_{t}(B)>0\}.
\end{equation}
\end{lemma}
\begin{pf}
Fix $\delta$ as in the lemma. Fix also $\omega$ such that $ X_{t}%
(B)>0$ and $ X_{s}(B)\rightarrow X_{t}(B)$ as
$s\uparrow t$. For this $\omega$, there is an $\varepsilon
\in(0,\delta)$ such that $X_{s}(B)>\varepsilon$ for all $s\in
(t-\varepsilon
,t)$. Hence%
%
%
\begin{equation}
\int_{t-\delta}^{t}ds\, X_{s}(B)f(t-s) \geq\varepsilon
\int_{t-\varepsilon}^{t}ds\, f(t-s) = \infty
\end{equation}
and we are done.
\end{pf}

Set
%
%
\begin{equation} \label{not.vy}%
\vartheta:= \frac{1}{1+\beta}
\end{equation}
and for $\varepsilon\in(0,t)$ let $\tau_{\varepsilon}(B)$ denote the first
moment in $(t-\varepsilon,t)$ in which a ``big
jump'' occurs. More precisely, define%
%
%
\begin{equation}
\tau_{\varepsilon}(B) := \inf\biggl\{ s\in(t-\varepsilon,t)\dvtx|\Delta
X_{s}|(B)>(t-s)^{\vartheta} \log^{\vartheta} \biggl(\frac{1}{t-s} \biggr) \biggr\}.
\end{equation}
\begin{lemma}[(Existence of big jumps)]\label{lem:3}
For $\varepsilon\in(0,t)$ and the open ball $B$,
%
%
\begin{equation} \label{equt:12}%
\mathbf{P} \bigl( \tau_{\varepsilon}(B)=\infty\bigr) \leq\mathbf{P}%
\bigl(X_{t}(B) = 0 \bigr).
\end{equation}
\end{lemma}
\begin{pf}
For simplicity, through the proof we write $\tau$ for $\tau
_{\varepsilon
}(B)$. It suffices to show that
%
%
\begin{equation}\label{equt:11}%
\mathbf{P} \{\tau=\infty, X_{t}(B)>0 \} = 0.
\end{equation}
To verify (\ref{equt:11}) we will mainly follow the lines of the proof of
Theorem 1.2(b) of~\cite{LeGallMytnik2003}. For
$ u\in(0,\varepsilon]$, define
\[
Z_{u}:=N \biggl((s,x,r)\dvtx s\in(t-\varepsilon, t-\varepsilon+u), x\in
B, r>(t-s)^{\vartheta}\log^{\vartheta} \biggl(\frac{1}{t-s} \biggr) \biggr)
\]
with the random measure $N$ introduced in Lemma \ref{L.mart.dec}(a). Then
%
%
\begin{equation} \label{equt:10}%
\{ \tau=\infty\} = \{Z_{\varepsilon}=0\}.
\end{equation}
Recall the formula for the compensator $\hat{N}$ of $N$ in
Lemma \ref{L.mart.dec}(b). From a classical time change result for
counting processes (see, e.g., Theorem 10.33 in \cite{Jakod1979}), we
get that
there exists a standard Poisson process $A=\{A(v)\dvtx v\geq0\}$ such
that
%
%
\begin{eqnarray}
Z_{u} & = & A \biggl(\varrho\int_{t-\varepsilon}^{t-\varepsilon+u}%
ds\, X_{s}(B)\int_{(t-s)^{\vartheta}\log^{\vartheta}
({1}/({t-s}))}^{\infty}dr\, r^{-2-\beta} \biggr)\nonumber\\[-8pt]\\[-8pt]
& = & A \biggl(\frac{\varrho}{1+\beta}\int_{t-\varepsilon}^{t-\varepsilon
+u}ds\, X_{s}(B) \frac{1}{(t-s)\log({1}/({t-s}))}\biggr),\nonumber
\end{eqnarray}
where we used notation (\ref{not.vy}). Then
%
%
\begin{eqnarray}\label{array}\quad
&& \mathbf{P}\bigl(Z_{\varepsilon}=0, X_{t}(B)>0\bigr) \nonumber\\[-8pt]\\[-8pt]
&&\qquad \leq\mathbf{P} \biggl(\int_{t-\varepsilon}^{t}ds\, X_{s}%
(B) \frac{1}{(t-s)\log({1}/({t-s}))}<\infty, X_{t}(B)>0 \biggr).\nonumber
\end{eqnarray}
It is easy to check that
%
%
\begin{equation}
\int_{t-\delta}^{t}ds\, \frac{1}{(t-s)\log({1}/({t-s}))}
= \infty\qquad\mbox{for all } \delta\in(0,\varepsilon).
\end{equation}
Therefore, by Lemma \ref{L.explosion},
%
%
\begin{equation}
\int_{t-\varepsilon}^{t}ds\, X_{s}(B) \frac{1}{(t-s)\log
({1}/({t-s}))} = \infty\qquad\mbox{on } \{X_{t}%
(B)>0\}.
\end{equation}
Thus, the probability in (\ref{array}) equals 0. Hence, together
with (\ref{equt:10}) claim (\ref{equt:11}) follows.
\end{pf}

Set $\varepsilon_{n}:=2^{-n}$, $n\geq1$. Then we choose
open balls $B_{n}\uparrow B$ such that $ $%
%
%
\begin{equation} \label{Bs}%
\overline{B_{n}}\subset B_{n+1}\subset B \quad\mbox{and}\quad \sup_{y\in
B^{\mathrm{c}}, x\in B_{n}, 0<s\leq\varepsilon_{n}} p_{s}^{\alpha
}(x-y) \mathop{\longrightarrow}_{n\uparrow\infty}0.
\end{equation}
Fix $n\geq1$ such that $ \varepsilon_{n}<t$. Define $\tau
_{n}:=\tau_{\varepsilon_{n}}(B_{n})$.

In order to get a lower bound for $ \Vert X_{t} \Vert_{B}$ we use
the following inequality:
%
%
\begin{equation} \label{123}%
\Vert X_{t} \Vert_{B} \geq\int_{B}dy\, X_{t}%
(y) p_{r}^{\alpha}(y-x),\qquad x\in B, r>0.
\end{equation}
On the event $\{\tau_{n}<t\}$, denote by $\zeta_{n}$ the spatial
location in
$B_{n}$ of the jump at time~$\tau_{n }$, and by $r_{n}$ the size
of the jump, meaning that $\Delta X_{\tau_{n}}=r_{n}\delta_{\zeta_{n}}%
$. Then specializing (\ref{123}),%
%
%
\begin{equation} \label{123'}%
\Vert X_{t} \Vert_{B} \geq\int_{B}dy\, X_{t}%
(y) p_{t-\tau_{n}}^{\alpha}(y-\zeta_{n}) \qquad\mbox{on the event } \{\tau
_{n}<t\}.
\end{equation}
From the strong Markov property at time $\tau_{n }$, together with the
branching property of superprocesses, we know that conditionally on
$\{\tau_{n}<t\}$, the process $\{X_{\tau_{n}+u}\dvtx u\geq 0\}$ is
bounded below in distribution by $\{\widetilde{X}_{u}^{n}\dvtx u\geq
0\}$ where $\widetilde{X}^{n}$ is a super-Brownian motion with initial
value $r_{n}\delta_{\zeta_{n}}$. Hence, from (\ref{123'}) we get
%
%
\begin{eqnarray} \label{firstbound}
&& \mathbf{E}\exp\{- \Vert X_{t} \Vert
_{B} \}\nonumber\\
&&\qquad \leq\mathbf{E} \mathsf{1}_{\{\tau_{n}<t\}}\exp\biggl\{ -\int
_{B}dy\, X_{t}(y) p_{t-\tau_{n}}^{\alpha}(y-\zeta_{n}) \biggr\}
+ \mathbf{P}(\tau_{n}=\infty)\nonumber\\[-8pt]\\[-8pt]
&&\qquad \leq\mathbf{E} \mathsf{1}_{\{\tau_{n}<t\}}\mathbf{E}_{r_{n}%
\delta_{\zeta_{n}}}\exp\biggl\{ -\int_{B}dy X_{t-\tau_{n}%
}(y) p_{t-\tau_{n}}^{\alpha}(y-\zeta_{n}) \biggr\}\nonumber\\
&&\qquad\quad{} + \mathbf{P}(\tau
_{n}=\infty).\nonumber
\end{eqnarray}
Note that on the event $\{\tau_{n}<t\}$, we have
%
%
\begin{equation} \label{not.hbeta}%
r_{n}\geq(t-\tau_{n})^{\vartheta}\log^{\vartheta} \biggl(\frac{1}{t-\tau_{n}%
} \biggr) =: h_{\beta}(t-\tau_{n}).
\end{equation}
We now claim that
%
%
\begin{equation} \label{secondlimit}\quad
\lim_{n\uparrow\infty} \sup_{0<s<\varepsilon_{n}, x\in B_{n}, r\geq
h_{\beta}(s)}\mathbf{E}_{r\delta_{x}}\exp\biggl\{ - \int_{B}%
dy\, X_{s}(y) p_{s}^{\alpha}(y-x) \biggr\} = 0.
\end{equation}
To verify (\ref{secondlimit}), let $ s\in(0,\varepsilon_{n})$,
$x\in B_{n}$ and $r\geq h_{\beta}(s)$. Then using the
Laplace transition functional of the superprocess we get
%
%
\begin{eqnarray} \label{Lapla}%
\mathbf{E}_{r\delta_{x}}\exp\biggl\{ -\int_{B}dy\, X_{s}%
(y) p_{s}^{\alpha}(y-x) \biggr\} &=& \exp\{ -r v_{s,x}%
^{n}(s,x) \} \nonumber\\[-8pt]\\[-8pt]
&\leq&\exp\{ -h_{\beta}(s) v_{s,x}^{n}(s,x) \},\nonumber
\end{eqnarray}
where the nonnegative function $ v_{s,x}^{n}=\{v_{s,x}^{n}(s^{\prime
},x^{\prime})\dvtx s^{\prime}>0$, $x^{\prime}\in\mathsf{R}^{d}\}$
solves the log-Laplace integral equation
%
%
\begin{eqnarray}\label{equt:4}\quad
v_{s,x}^{n}(s^{\prime}, x^{\prime}) &=& \int_{\mathsf{R}^{d}}%
dy\, p_{s^{\prime}}^{\alpha}(y-x^{\prime}) 1_{B}(y) p_{s}^{\alpha
}(y-x)\nonumber\\
&&{} +\int_{0}^{s^{\prime}}dr^{\prime}\int_{\mathsf{R}^{d}}%
dy\, p_{s^{\prime}-r^{\prime}}^{\alpha}(y-x^{\prime}) [av_{s,x}%
^{n}(r^{\prime},y)\\
&&\hspace*{141.3pt}{}-b(v_{s,x}^{n}(r^{\prime},y))^{1+\beta
} ]\nonumber
\end{eqnarray}
related to (\ref{logLap}).
\begin{lemma}[(Another explosion)]\label{lem:2} Under the conditions $d>1$
or $\alpha\leq1+\beta$, we have
%
%
\begin{equation} \label{limitv}%
\lim_{n\uparrow\infty} \Bigl(\inf_{0<s<\varepsilon_{n}, x\in B_{n}}h_{\beta
}(s) v_{s,x}^{n}(s,x) \Bigr) = +\infty.
\end{equation}
\end{lemma}

Let us postpone the proof of Lemma \ref{lem:2}.
\begin{pf*}{Completion of Proof of Theorem
\ref{T.prop.dens}\textup{(c)}}
Our claim (\ref{secondlimit}) readily follows from
estimate (\ref{Lapla}) and (\ref{limitv}). Moreover, according to
(\ref{secondlimit}), by passing to the limit $n\uparrow\infty$ in the
right-hand side of (\ref{firstbound}), and then using Lemma \ref
{lem:3}, we arrive
at
%
%
\begin{equation}\qquad
\mathbf{E}\exp\{- \Vert X_{t} \Vert
_{B} \} \leq\limsup_{n\uparrow\infty}\mathbf{P} ( \tau
_{n}=\infty) \leq\limsup_{n\uparrow\infty}\mathbf{P} \bigl(
X_{t}(B_{n})=0 \bigr).
\end{equation}
Since the event $\{X_{t}(B)=0\}$ is the nonincreasing limit as
$n\uparrow\infty$ of the events $\{X_{t}(B_{n})=0\}$ we get%
%
%
\begin{equation}
\mathbf{E}\exp\{- \Vert X_{t} \Vert
_{B} \} \leq\mathbf{P} \bigl(X_{t}%
(B) = 0 \bigr).
\end{equation}
Since obviously $ \Vert X_{t} \Vert_{B}=0$ if and only if
$X_{t}(B)=0$, we see that (\ref{basicprop}) follows from this last
bound. The
proof of Theorem 1(c) is finished for $U=B$.
\end{pf*}
\begin{pf*}{Proof of Lemma \ref{lem:2}}
We start with a determination of the asymptotics of the first term at
the right-hand side of the log-Laplace equation (\ref{equt:4}) at
$(s^{\prime},x^{\prime})=(s,x)$. Note that
%
%
\begin{eqnarray}\label{equt:i1}\quad
&& \int_{\mathsf{R}^{d}}dy\, p_{s}^{\alpha}(y-x) 1_{B}(y) p_{s}%
^{\alpha}(y-x)\nonumber\\[-8pt]\\[-8pt]
&&\qquad = \int_{\mathsf{R}^{d}}dy\, p_{s}^{\alpha}(y-x) p_{s}%
^{\alpha}(y-x)-\int_{B^{\mathrm{c}}}dy\, p_{s}^{\alpha}(y-x) p_{s}%
^{\alpha}(y-x).\nonumber
\end{eqnarray}
In the latter formula line, the first term equals $ p_{2s}^{\alpha
}(0)=Cs^{-d/\alpha}$ whereas the second one is bounded from above
by%
%
%
\begin{equation} \label{equt:i2}%
\sup_{0<s<\varepsilon_{n}, x\in B_{n}, y\in B^{\mathrm{c}}}p_{s}^{\alpha
}(y-x) \mathop{\longrightarrow}_{n\uparrow\infty}0,
\end{equation}
where the last convergence follows by assumption (\ref{Bs}) on $B_{n }%
$. Hence from (\ref{equt:i1}) and (\ref{equt:i2}) we obtain
%
%
\begin{equation}\label{equt:i3}\quad
\int_{\mathsf{R}^{d}}dy\, p_{s}^{\alpha}(y-x) 1_{B}(y) p_{s}%
^{\alpha}(y-x) = C s^{-d/\alpha}+\mathrm{o}(1) \qquad\mbox{as } n\uparrow
\infty,
\end{equation}
uniformly in $ s\in(0,\varepsilon_{n})$ and $ x\in B_{n }$.

To simplify notation, we write $v^{n}:=v_{s,x}^{n}$. Next,
from (\ref{equt:4}) we can easily get the upper bound%
%
%
\begin{eqnarray}\label{equt:i4}
v^{n}(s^{\prime},x^{\prime}) &\leq& e^{|a|s^{\prime}}\int
_{\mathsf{R}^{d}}dy\, p_{s^{\prime}}^{\alpha}(y-x^{\prime}%
) p_{s}^{\alpha}(y-x) \nonumber\\[-8pt]\\[-8pt]
&=& e^{|a|s^{\prime}} p_{s^{\prime}%
+s}^{\alpha}(x-x^{\prime}).\nonumber
\end{eqnarray}
Then we have
%
%
\begin{eqnarray}\label{equt:i5}\quad
&& \int_{0}^{s}dr^{\prime}\int_{\mathsf{R}^{d}}d%
y\, p_{s-r^{\prime}}^{\alpha}(y-x) (v^{n}(r^{\prime},y))^{1+\beta}\nonumber\\
&&\qquad \leq e^{|a|(1+\beta)s} \int_{0}^{s}dr^{\prime}%
\int_{\mathsf{R}^{d}}dy\, p_{s-r^{\prime}}^{\alpha}%
(y-x)
\bigl(p_{r^{\prime}+s}^{\alpha}(x-y)\bigr)^{1+\beta}\nonumber\\[-8pt]\\[-8pt]
&&\qquad \leq e^{|a|(1+\beta)s}(p_{s}^{\alpha}(0))^{\beta
}\int_{0}^{s}dr^{\prime}\int_{\mathsf{R}^{d}}d%
y\, p_{s-r^{\prime}}^{\alpha}(y-x) p_{r^{\prime}+s}^{\alpha}(x-y)\nonumber
\nonumber\\
&&\qquad = e^{|a|(1+\beta)s}(p_{s}^{\alpha}(0))^{\beta}%
\int_{0}^{s}dr^{\prime}\, p_{2s}^{\alpha}(0) = C e%
^{|a|(1+\beta)s}s^{1-d(1+\beta)/\alpha}\nonumber
\end{eqnarray}
and, similarly,%
%
%
\begin{equation} \label{added}%
\int_{0}^{s}dr^{\prime}\int_{\mathsf{R}^{d}}d%
y\, p_{s-r^{\prime}}^{\alpha}(y-x) av^{n}(r^{\prime},y) \geq
-C |a| e^{|a|s} s^{1-d/\alpha}.
\end{equation}
Summarizing, by (\ref{equt:4}), (\ref{equt:i3}), (\ref{equt:i5}) and
(\ref{added}),%
%
%
\begin{equation}\label{equt:i6}\hspace*{28pt}
v^{n}(s,x)\geq C s^{-d/\alpha}+\mathrm{o}(1)-C e^{|a|(1+\beta
)s} s^{1-d(1+\beta)/\alpha}-C |a| e^{|a|s} s^{1-d/\alpha}
\end{equation}
uniformly in $ s\in(0,\varepsilon_{n})$ and $ x\in B_{n }$.
According to our general assumption $d<\alpha/\beta$, we conclude that the
right-hand side of (\ref{equt:i6}) behaves like $Cs^{-d/\alpha}$ as
$s\downarrow0$ uniformly in $ s\in(0,\varepsilon_{n})$. Now recalling
definitions (\ref{not.hbeta}) and (\ref{not.vy}) as well as our assumption
that $d>1$ or $\alpha\leq1+\beta$, we immediately get
%
%
\begin{equation} \label{equt:i9}%
\lim_{n\uparrow\infty} \inf_{0<s<\varepsilon_{n}}h_{\beta}(s)
s^{-d/\alpha
} = +\infty.
\end{equation}
By (\ref{equt:i6}), this implies (\ref{limitv}), and the proof of the
lemma is complete.
\end{pf*}

\section[Optimal local H\"{o}lder index: Proof of Theorem 1.2(b)]{Optimal local
H\"{o}lder index: Proof of Theorem \protect\ref{T.prop.dens}\textup{(b)}}\label{sec:5}

We return to \mbox{$d=1$} and continue to assume that $ t>0$ and
$\mu\in\mathcal{M}_{\mathrm{f}}\setminus\{0\}$ are fixed. In the proof
of Theorem \ref{T.prop.dens}(b) we implement the following idea. We
show that there exists a sequence of ``big'' jumps of $X$ that occur
close to time $t$ and these jumps in fact destroy the local H\"{o}lder
continuity of any index greater or equal than $\eta_{\mathrm{c} }$.

As in the proof of Theorem \ref{T.prop.dens}(c) in the previous
section, we may work with a fixed open interval $U$. For simplicity we
consider $ U=(0,1)$. Put
%
%
\begin{equation}
I_{k}^{(n)} := \biggl[\frac{k}{2^{n}},\frac{k+1}{2^{n}} \biggr),\qquad
n\geq1, 0\leq k\leq2^{n}-1.
\end{equation}
Choose $n_{0}$ such that $2^{-\alpha n_{0}}<t$. For $ n\geq n_{0}%
$ and $ 2\leq k\leq2^{n}+1$, denote by $A_{n,k}$ the
following event:
%
%
\begin{equation} \label{c5.2}\hspace*{28pt}
\biggl\{ \Delta X_{s}\bigl(I_{k-2}^{(n)}\bigr) \geq\frac{c_{\mazinti{(\ref{c5.2})}}}%
{2^{{\alpha}/({1+\beta})n}} n^{{1}/({1+\beta})} \mbox{ for some }%
s\in\bigl[t-2^{-\alpha n}, t-2^{-\alpha(n+1)}\bigr) \biggr\}
\end{equation}
with $ c_{\mazinti{(\ref{c5.2})}}:=(\alpha2^{-\alpha}\log2)^{{1}/({1+\beta})}$,
and for $N\geq n_{0}$ write
%
%
\begin{equation}
\widetilde{A}_{N} := \bigcup_{n=N}^{\infty}\bigcup
_{k=2}^{2^{n}+1}A_{n,k }.
\end{equation}
\begin{lemma}[(Again existence of big jumps)]\label{lem:1'}
For any $N\geq n_{0 }$,
%
%
\begin{equation}
\mathbf{P}\{\widetilde{A}_{N} | X_{t}(U)>0\} = 1.
\end{equation}
\end{lemma}
\begin{pf}
For $ s\in[t-2^{-\alpha n}, t-2^{-\alpha(n+1)})$ we have%
%
%
\begin{eqnarray}
\biggl((t-s)\log\biggl(\frac{1}{t-s}\biggr) \biggr)^{{1}/({1+\beta})}
&\geq&\bigl(2^{-\alpha(n+1)}\log2^{\alpha n}
\bigr)^{{1}/({1+\beta})}\nonumber\\[-8pt]\\[-8pt]
&=& c_{\mazinti{(\ref{c5.2})}} 2^{-{\alpha}/({1+\beta}) n}
n^{{1}/({1+\beta})}.\nonumber
\end{eqnarray}
Therefore,%
\begin{eqnarray*}
\bigcup_{k=2}^{2^{n}+1} A_{n,k}&\supseteq&\biggl\{\Delta X_{s}(U) \geq
\biggl((t-s)\log\biggl(\frac{1}{t-s}\biggr) \biggr)^{{1}/({1+\beta})} \\
&&\hspace*{17.6pt}\mbox{for some }s\in\bigl[t-2^{-\alpha n}, t-2^{-\alpha(n+1)}\bigr)
\biggr\}
\end{eqnarray*}
and, consequently,%
%
%
\begin{eqnarray}\hspace*{28pt}
\widetilde{A}_{N} &=& \bigcup_{n=N}^{\infty}\bigcup
_{k=2}^{2^{n}+1}A_{n,k}\nonumber\\[-8pt]\\[-8pt]
&\supseteq&\biggl\{ \Delta X_{s}(U)\geq\biggl((t-s)\log\biggl(\frac{1}%
{t-s}\biggr) \biggr)^{{1}/({1+\beta})} \mbox{ for some } s\geq
t-2^{-N} \biggr\}\nonumber
\end{eqnarray}
and we are done by Lemma \ref{lem:3}.
\end{pf}

Now we are going to define increments of $Z_{t}^{2}$ on the dyadic sets
$\{\frac{k}{2^{n}}\dvtx k=0,\ldots,2^{n}\}$. By definition
(\ref{rep.dens}),
%
%
\begin{eqnarray}\quad
&& Z_{t}^{2} \biggl(\frac{k}{2^{n}} \biggr) -Z_{t}^{2}
\biggl(\frac{k+1}{2^{n}} \biggr) \nonumber\\
&&\qquad =\int_{(0,t]\times\mathsf{R}} M (d%
(s,y) ) \biggl(p_{t-s}^{\alpha} \biggl(\frac{k}{2^{n}}-y \biggr) -p_{t-s}^{\alpha}
\biggl(\frac{k+1}{2^{n}}-y \biggr)
\biggr) \nonumber\\[-8pt]\\[-8pt]
&&\qquad =\int_{(0,t]\times\mathsf{R}} M (d%
(s,y) ) \biggl(p_{t-s}^{\alpha}
\biggl(\frac{k}{2^{n}}-y \biggr) -p_{t-s}^{\alpha}
\biggl(\frac{k+1}{2^{n}}-y \biggr)
\biggr) _{+}\nonumber\\
&&\qquad\quad{} +\int_{(0,t]\times\mathsf{R}} M
(d(s,y) ) \biggl(p_{t-s}^{\alpha} \biggl(
\frac{k}{2^{n}}-y \biggr) -p_{t-s}^{\alpha}
\biggl(\frac{k+1}{2^{n}}-y \biggr) \biggr) _{-}.\nonumber
\end{eqnarray}
Then according to Lemma \ref{L9} there exist spectrally positive stable
processes $L_{n,k}^{+}$ and $L_{n,k}^{-}$ of index $1+\beta$ such that
%
%
\begin{equation}
Z_{t}^{2} \biggl(\frac{k}{2^{n}} \biggr) -Z_{t}^{2} \biggl(\frac{k+1}{2^{n}%
} \biggr) = L_{n,k}^{+}(T_{+})-L_{n,k}^{-}(T_{-}),
\end{equation}
where
%
%
\begin{equation} \label{not.Tpm}%
T_{\pm} :=\int_{0}^{t}ds\int_{\mathsf{R}}X_{s}(dy)
\biggl(p_{t-s}^{\alpha} \biggl(\frac{k}{2^{n}}-y \biggr)
-p_{t-s}^{\alpha} \biggl(\frac{k+1}{2^{n}}-y \biggr) \biggr)
_{\pm}^{ 1+\beta}.
\end{equation}
Fix $\varepsilon\in(0,\frac{1}{1+\beta})$ for a while.
Let us define the following events:
%
%
\begin{eqnarray} \label{5.10}
B_{n,k} &:=& \bigl\{ L_{n,k}^{+}(T_{+})\geq2^{-\eta_{\mathrm{c}}n}n^{
{1}/({1+\beta})-\varepsilon} \bigr\} \cap\{ L_{n,k}^{-}(T_{-})\leq
2^{-\eta_{\mathrm{c}}n-\varepsilon n} \}\nonumber\\[-8pt]\\[-8pt]
&=:& B_{n,k}^{+}\cap B_{n,k }^{-}\nonumber
\end{eqnarray}
(with notation in the obvious correspondence). Define the following event:
%
%
\begin{eqnarray}
D_{N} :\!&=& \bigcup_{n=N}^{\infty}\bigcup_{k=2}^{2^{n}+1} (
A_{n,k}\cap B_{n,k} ) \nonumber\\[-8pt]\\[-8pt]
&\supseteq& \bigcup_{n=N}^{\infty}\bigcup_{k=2}^{2^{n}+1}A_{n,k}%
\biggm\backslash\bigcup_{n=N}^{\infty}\bigcup_{k=2}^{2^{n}+1} ( A_{n,k}\cap
B_{n,k}^{\mathrm{c}} ) .\nonumber
\end{eqnarray}
An estimation of the probability of $D_{N}$ is crucial for the proof of
Theorem \ref{T.prop.dens}(b). In fact we are going to show that conditionally
on $\{X_{t}(U)>0\}$, the event $D_{N}$ happens with probability one for any
$N$. This in turn implies that for any $N$ one can find $n\geq N$ sufficiently
large such that there exists an interval $[\frac{k}{2^{n}},\frac{k+1}{2^{n}}]$
on which the increment $Z_{t}^{2}(\frac{k}{2^{n}})-Z_{t}^{2}(\frac
{k+1}{2^{n}%
})$ is of order $L_{n,k}^{+}(T_{+})\geq2^{-\eta_{\mathrm{c}}n}n^{
{1}/({1+\beta})-\varepsilon}$ [since the other term $L_{n,k}^{-}(T_{-})$
is much smaller on that interval]. This implies the statement of
Theorem \ref{T.prop.dens}(b). Detailed arguments follow.

By Lemma \ref{lem:1'} we get
%
%
\begin{equation} \label{13}\quad
\mathbf{P} \{ D_{N} | X_{t}(U)>0 \} \geq1-\mathbf{P}%
\Biggl\{\bigcup_{n=N}^{\infty}\bigcup_{k=2}^{2^{n}+1} ( A_{n,k}\cap
B_{n,k}^{\mathrm{c}} ) \bigg| X_{t}(U)>0 \Biggr\}.
\end{equation}
Recall $A^{\varepsilon}$ defined in (\ref{A.eps}). Note that
%
%
\begin{eqnarray} \label{16}
&& \mathbf{P} \Biggl( \bigcup_{n=N}^{\infty}\bigcup_{k=2}^{2^{n}+1}
(A_{n,k}\cap B_{n,k}^{\mathrm{c}} ) \Biggr)\nonumber\\
&&\qquad \leq\mathbf{P}(A^{\varepsilon,\mathrm{c}})+\mathbf{P}
\Biggl(\bigcup_{n=N}^{\infty}\bigcup_{k=2}^{2^{n}+1} ( A^{\varepsilon}\cap
A_{n,k}\cap B_{n,k}^{\mathrm{c}} ) \Biggr) \\
&&\qquad \leq2\varepsilon+\mathbf{P} \Biggl( \bigcup_{n=N}^{\infty}\bigcup
_{k=2}^{2^{n}+1} ( A^{\varepsilon}\cap A_{n,k}\cap B_{n,k}^{\mathrm{c}%
} ) \Biggr) .\nonumber
\end{eqnarray}
\begin{lemma}[(Probability of small increments)]\label{Prop2}
For all $ \varepsilon>0$ sufficiently small,%
%
%
\begin{equation} \label{5.1}%
\lim_{N\uparrow\infty}\mathbf{P} \Biggl( \bigcup_{n=N}^{\infty}\bigcup
_{k=2}^{2^{n}+1} ( A^{\varepsilon}\cap A_{n,k}\cap B_{n,k}^{\mathrm{c}%
} ) \Biggr) = 0.
\end{equation}
\end{lemma}

We postpone the proof of this lemma to the end of this section. Instead we
will show now, how it implies Theorem \ref{T.prop.dens}(b).
\begin{pf*}{Completion of proof of Theorem
\ref{T.prop.dens}\textup{(b)}}
From Lemma \ref{Prop2} and (\ref{16}) it follows that%
%
%
\begin{equation}\hspace*{4pt}
\limsup_{N\uparrow\infty}\mathbf{P} \Biggl\{\bigcup_{n=N}^{\infty}\bigcup
_{k=2}^{2^{n}+1} ( A_{n,k}\cap B_{n,k}^{\mathrm{c}} )
\bigg| X_{t}(U)>0 \Biggr\} \leq\frac{2\varepsilon}{\mathbf{P}%
(X_{t}(U)>0)} .
\end{equation}
Since $\varepsilon$ can be arbitrarily small, the latter $\limsup$ expression
equals 0. Combining this with estimate (\ref{13}), we get%
%
%
\begin{equation}
\lim_{N\uparrow\infty}\mathbf{P} \{ D_{N} | X_{t}(U)>0 \}
= 1.
\end{equation}
Since $ D_{N}\downarrow\bigcap_{N=n_{0}}^{\infty}D_{N}=:D_{\infty}$
as $N\uparrow\infty$, we conclude that
%
%
\begin{equation}
\mathbf{P} \{ D_{\infty} | X_{t}(U)>0 \} = 1.
\end{equation}
This means that, almost surely on $\{X_{t}(U)>0\}$, there is a
sequence $(n_{j},k_{j})$ such that
%
%
\begin{equation}
Z_{t}^{2}\biggl(\frac{k_{j}}{2^{n_{j}}}\biggr)-Z_{t}^{2}\biggl(\frac{k_{j}%
+1}{2^{n_{j}}}\biggr) \geq2^{-\eta_{\mathrm{c}} n_{j}} n_{j}^{
{1}/({1+\beta})-\varepsilon}.
\end{equation}
This inequality implies the claim in Theorem \ref{T.prop.dens}(b).
\end{pf*}

We now prepare for the proof of Lemma \ref{Prop2}.
Actually by using (\ref{5.10}), we represent the probability in (\ref{5.1})
as a sum of the two following probabilities:
%
%
\begin{eqnarray} \label{5.2}%
&& \mathbf{P} \Biggl( \bigcup_{n=N}^{\infty}\bigcup_{k=2}^{2^{n}+1}
(A^{\varepsilon}\cap A_{n,k}\cap B_{n,k}^{\mathrm{c}})\Biggr)
\nonumber\\
&&\qquad = \mathbf{P} \Biggl( \bigcup_{n=N}^{\infty}\bigcup_{k=2}^{2^{n}%
+1} ( A^{\varepsilon}\cap A_{n,k}\cap B_{n,k}^{+,\mathrm{c}} )
\Biggr) \\
&&\qquad\quad{} + \mathbf{P} \Biggl( \bigcup_{n=N}^{\infty}\bigcup_{k=2}%
^{2^{n}+1} ( A^{\varepsilon}\cap A_{n,k}\cap B_{n,k}^{-,\mathrm{c}%
} ) \Biggr) .\nonumber
\end{eqnarray}
Now we will handle each term on the right-hand side of (\ref{5.2}) separately.
\begin{lemma}[{[First term in (\ref{5.2})]}]\label{L.1part}
For $\varepsilon\in(0,\frac{1}{1+\beta})$,%
%
%
\begin{equation}
\lim_{N\uparrow\infty}\mathbf{P} \Biggl( \bigcup_{n=N}^{\infty}\bigcup
_{k=2}^{2^{n}+1} ( A^{\varepsilon}\cap A_{n,k}\cap B_{n,k}^{+,\mathrm{c}
} ) \Biggr) = 0.
\end{equation}
\end{lemma}
\begin{pf}
Consider the process $L_{n,k}^{+}(s), s\leq T_{+}$. On $A_{n,k}$ there
exists a jump of the martingale measure $M$ of the form $ r^{\ast}%
\delta_{s^{\ast},y^{\ast}}$ for some
%
%
\begin{eqnarray}
r^{\ast}&\geq& c_{\mazinti{(\ref{c5.2})}} 2^{-{\alpha}/({1+\beta})n}n^{{1}
/({1+\beta})},\nonumber\\[-8pt]\\[-8pt]
s^{\ast}&\in&\bigl[ t-2^{-\alpha n}, t-2^{-\alpha
(n+1)}\bigr],\qquad y^{\ast}\in I_{k-2 }^{(n)}.\nonumber
\end{eqnarray}
Hence
%
%
\begin{eqnarray}\label{equt:2}\hspace*{32pt}
\Delta L_{n,k}^{+}(s^{\ast})&\geq&\inf_{y\in I_{k-2}^{(n)}, s\in[
2^{-\alpha(n+1)},2^{-\alpha n}]}  \biggl(p_{s}%
^{\alpha} \biggl(\frac{k}{2^{n}}-y \biggr) -p_{s}^{\alpha}
\biggl(\frac{k+1}{2^{n}}-y \biggr) \biggr) _{+}\nonumber\\[-8pt]\\[-8pt]
&&\hspace*{93pt}{} \times c_{\mazinti{(\ref{c5.2})}} 2^{-{\alpha}/({1+\beta})n} n^{
{1}/({1+\beta})}.\nonumber
\end{eqnarray}
It is easy to get
%
%
\begin{eqnarray}\label{equt:1'}\qquad
&& \inf_{y\in I_{k-2}^{(n)}, s\in[2^{-\alpha(n+1)},2^{-\alpha n}%
]} \biggl(p_{s}^{\alpha} \biggl(\frac{k}{2^{n}%
}-y \biggr) -p_{s}^{\alpha} \biggl(\frac{k+1}{2^{n}}-y \biggr) \biggr)
_{+}\nonumber\\
&&\qquad = \mathop{\inf_{
2^{-n}\leq z\leq2^{-n+1},}}_
{s\in[2^{-\alpha(n+1)},2^{-\alpha n}]}
\bigl(p_{s}^{\alpha} ( z ) -p_{s}^{\alpha
} ( z+2^{-n} ) \bigr) _{+}\nonumber\\
&&\qquad = \mathop{\inf_{2^{-n}\leq z\leq2^{-n+1},}}_
{s\in[2^{-\alpha(n+1)},2^{-\alpha n}]}
s^{-1/\alpha} \bigl( p_{1}^{\alpha}(zs^{-1/\alpha})-p_{1}^{\alpha
}\bigl((z+2^{-n})s^{-1/\alpha}\bigr) \bigr) _{+}\\
&&\qquad \geq2^{n}{\mathop{\inf_{
2^{-n}\leq z\leq3\cdot2^{-n},}}_{s\in[2^{-\alpha(n+1)},2^{-\alpha n}]}}
\vert(p_{1}^{\alpha})^{\prime}(zs^{-1/\alpha}) \vert
2^{-n}s^{-1/\alpha}\nonumber\\
&&\qquad \geq{2^{n}\inf_{1\leq x\leq6} }\vert(p_{1}^{\alpha
})^{\prime} ( x ) \vert=: c_{\mazinti{(\ref{equt:1'})}}2^{n},\nonumber
\end{eqnarray}
where $c_{\mazinti{(\ref{equt:1'})}}>0$. In fact, from (\ref{L1.2}),%
%
%
\begin{equation}
\frac{d}{dz}p_{1}^{\alpha}(z) = -\int_{0}^{\infty
}ds\, q_{1}^{\alpha/2}(s) \frac{z}{2s} p_{s}^{(2)}(z) \neq0,\qquad
z\neq0,
\end{equation}
and $(p_{\alpha}^{(2)})^{\prime}(x)\not=0$ for any $x\not=0$.
Apply (\ref{equt:1'}) in (\ref{equt:2}) to arrive at
%
%
\begin{equation} \label{equt:3}%
\Delta L_{n,k}^{+}(s^{\ast}) \geq c_{\mazinti{(\ref{equt:3})}}2^{(1-{\alpha
}/({1+\beta}))n }n^{{1}/({1+\beta})} = c_{\mazinti{(\ref{equt:3})}}2^{-\eta
_{\mathrm{c}}n} n^{{1}/({1+\beta})}.
\end{equation}
Using Lemma \ref{L7} with $ \theta=1+\beta$ and%
%
%
\begin{equation}
\delta= (1+\beta)\mathsf{1}_{\{2\beta<\alpha-1\}}+(\alpha-\beta
-\varepsilon)\mathsf{1}_{\{2\beta\geq\alpha-1\} },
\end{equation}
we get, with $c_{\varepsilon}$ appearing in definition (\ref{A.eps}) of
$A^{\varepsilon}$,
%
%
\begin{equation} \label{not.tn}\hspace*{28pt}
T_{\pm} \leq c_{\varepsilon} \bigl(2^{-n(1+\beta)}\mathsf{1}_{\{2\beta
<\alpha-1\}}+2^{-n(\alpha-\beta-\varepsilon)}\mathsf{1}_{\{2\beta\geq
\alpha-1\}} \bigr) =: t_{n} \qquad\mbox{on } A^{\varepsilon}.
\end{equation}
Hence for all $n$ sufficiently large we obtain
%
%
\begin{eqnarray}\label{equt:4'}\quad
&& \mathbf{P} \bigl( L_{n,k}^{+}(T_{+})<2^{-\eta_{\mathrm{c}}n} n^{{1}/({1+\beta})
-\varepsilon}, A^{\varepsilon}\cap A_{n,k} \bigr) \nonumber\\
&&\qquad \leq\mathbf{P} \bigl( L_{n,k}^{+}(T_{+})<2^{-\eta_{\mathrm{c}}%
n} n^{{1}/({1+\beta})-\varepsilon},\nonumber\\
&&\hspace*{43.3pt} \Delta L_{n,k}^{+}(s^{\ast})\geq
c_{\mazinti{(\ref{equt:3})}}2^{-\eta_{\mathrm{c}}n} n^{{1}/({1+\beta})},
A^{\varepsilon} \bigr) \nonumber\\
&&\qquad \leq\mathbf{P} \biggl( \inf_{s\leq T^{+}}L_{n,k}^{+}(s)<-\frac{1}%
{2} c_{\mazinti{(\ref{equt:3})}}2^{-\eta_{\mathrm{c}}n} n^{{1}/({1+\beta})%
}, A^{\varepsilon} \biggr) \nonumber\\[-8pt]\\[-8pt]
&&\qquad \leq\mathbf{P} \biggl( \inf_{s\leq t_{n}}L_{n,k}^{+}(s)<-\frac{1}%
{2} c_{\mazinti{(\ref{equt:3})}}2^{-\eta_{\mathrm{c}}n} n^{{1}/({1+\beta})} \biggr)
\nonumber\\
&&\qquad \leq\exp\bigl\{-c_{\beta}(t_{n})^{-1/\beta}\bigl(c_{\mazinti{(\ref{equt:3})}}%
2^{-\eta_{\mathrm{c}}n} n^{{1}/({1+\beta})}\bigr)^{(1+\beta)/\beta
} \bigr\}\nonumber\\
&&\qquad \leq\exp\bigl\{-c_{\varepsilon}n^{1/\beta}\bigl(t_{n}^{-1}2^{-\eta_{\mathrm{c}
}(1+\beta)n}\bigr)^{1/\beta} \bigr\} \nonumber\\
&&\qquad\leq\exp\bigl\{-c_{\varepsilon}n^{1/\beta
}2^{(1-\varepsilon)n}\bigr\},\nonumber
\end{eqnarray}
where (\ref{equt:4'}) follows by Lemma \ref{L.small.values}, and the
rest is
simple algebra. From this we get that for $N$ sufficiently large,
%
%
\begin{eqnarray}\qquad
\mathbf{P} \Biggl( \bigcup_{n=N}^{\infty}\bigcup_{k=2}^{2^{n}+1}
(A^{\varepsilon}\cap A_{n,k}\cap B_{n,k}^{+,\mathrm{c}} ) \Biggr)
&\leq&
\sum_{n=N}^{\infty}\sum_{k=2}^{2^{n}+1}\mathbf{P} \bigl( A^{\varepsilon
}\cap A_{n,k}\cap B_{n,k}^{+,\mathrm{c}} \bigr) \nonumber\\
&\leq&
\sum_{n=N}^{\infty}\sum_{k=2}^{2^{n}+1}\exp\bigl\{-c_{\varepsilon
}n^{1/\beta} 2^{(1-\varepsilon)n}\bigr\} \\
&=& \sum_{n=N}^{\infty}2^{n}%
\exp\bigl\{-c_{\varepsilon}n^{1/\beta} 2^{(1-\varepsilon)n}\bigr\},\nonumber
\end{eqnarray}
which converges to $0$ as $N\uparrow\infty$, and we are done with the
proof of
Lemma \ref{L.1part}.
\end{pf}
\begin{lemma}[{[Second term in (\ref{5.2})]}]\label{prop:2}
For all $ \varepsilon>0$ sufficiently small,
%
%
\begin{equation}
\lim_{N\uparrow\infty}\mathbf{P} \Biggl( \bigcup_{n=N}^{\infty}\bigcup
_{k=2}^{2^{n}+1} ( A^{\varepsilon}\cap A_{n,k}\cap
B_{n,k}^{-,\mathrm{c}}) \Biggr) = 0.
\end{equation}
\end{lemma}

The proof of this lemma will be postponed almost to the end of the section.
For its preparation, fix $\rho\in(0,\frac{1}{2})$. Define
\[
A_{n}^{\rho} := \Bigl\{\omega\dvtx\mbox{there exists }I_{k}^{(n)} \mbox{ with
}\sup_{s\in[ t-2^{-\alpha(1-\rho)n}, t)}X_{s}\bigl(I_{k}^{(n)}%
\bigr)\geq2^{-n(1-2\rho)} \Bigr\}.
\]
Note that%
%
%
\begin{eqnarray}
&& \mathbf{P} \Biggl( \bigcup_{k=2}^{2^{n}+1} ( A^{\varepsilon}\cap
A_{n,k}\cap B_{n,k}^{-,\mathrm{c}} ) \Biggr) \nonumber\\
&&\qquad \leq\mathbf{P}(A_{n}^{\rho})+\mathbf{P} \Biggl( \bigcup_{k=2}^{2^{n}%
+1} ( A_{n}^{\rho,\mathrm{c}}\cap A^{\varepsilon}\cap A_{n,k}\cap
B_{n,k}^{-,\mathrm{c}} ) \Biggr) \\
&&\qquad \leq\mathbf{P}(A_{n}^{\rho})+\sum_{k=2}^{2^{n}+1}\mathbf{P} (
A_{n}^{\rho,\mathrm{c}}\cap A^{\varepsilon}\cap A_{n,k}\cap B_{n,k}%
^{-,\mathrm{c}} ).\nonumber
\end{eqnarray}
Now let us introduce the notation
%
%
\begin{equation}
B_{n,k}^{-,1} := \Bigl\{\sup_{s\leq T_{-}}\Delta L_{n,k}^{-}(s)\leq
2^{-\eta_{\mathrm{c}}n-\varepsilon n} \Bigr\}.
\end{equation}
Then we have
%
%
\begin{eqnarray} \label{3terms}%
&& \mathbf{P} \Biggl( \bigcup_{k=2}^{2^{n}+1} ( A^{\varepsilon}\cap
A_{n,k}\cap B_{n,k}^{-,\mathrm{c}} ) \Biggr) \nonumber\\
&&\qquad \leq\mathbf{P}(A_{n}^{\rho})+\sum_{k=2}^{2^{n}+1}\mathbf{P} (
A_{n}^{\rho,\mathrm{c}}\cap A^{\varepsilon}\cap A_{n,k}\cap B_{n,k}%
^{-,\mathrm{c}} ) \nonumber\\
&&\qquad \leq\mathbf{P}(A_{n}^{\rho})+\sum_{k=2}^{2^{n}+1}\mathbf{P} (
A^{\varepsilon}\cap B_{n,k}^{-,\mathrm{c}}\cap B_{n,k}^{-,1} )
\\
&&\qquad\quad{} +\sum_{k=2}^{2^{n}+1}\mathbf{P} ( A^{\varepsilon}\cap
A_{n}^{\rho,\mathrm{c}}\cap A_{n,k}\cap B_{n,k}^{-,1,\mathrm{c}} )
\nonumber\\
&&\qquad =: \mathbf{P}(A_{n}^{\rho})+\sum_{k=2}^{2^{n}+1}P_{n,k}^{\varepsilon}%
+\sum_{k=2}^{2^{n}+1}P_{n,k }^{\varepsilon,\varrho}.\nonumber
\end{eqnarray}
In the following lemmas we consider the three terms in (\ref{3terms})
separately.
\begin{lemma}[{[First term in (\ref{3terms})]}]\label{L.small.prob}
There exists a constant $c_{\mazinti{(\ref{small.prob})}}$ independent
of $\rho\in(0,\frac{1}{2})$ such that
%
%
\begin{equation} \label{small.prob}%
\mathbf{P}(A_{n}^{\rho}) \leq c_{\mazinti{(\ref{small.prob})}} 2^{-\rho
n},\qquad n\geq n_{0 }.
\end{equation}
\end{lemma}
\begin{pf}
Fix $n\geq n_{0 }$. Define the stopping time $ \tau_{n}=\tau_{n}(\rho
)$ as%
%
%
\begin{equation}
\inf\bigl\{s\in\bigl[ t-2^{-\alpha(1-\rho)n},t\bigr)\dvtx X_{s}\bigl(I_{k}^{(n)}\bigr)
\geq2^{-n(1-2\rho)}\mbox{ for some }I_{k}^{(n)} \bigr\},
\end{equation}
if $ \omega\in A_{n}^{\rho}$, and as $t$ if $ \omega\in
A_{n }^{\rho,\mathrm{c}}$. Fix any $\omega\in A_{n}^{\rho}%
$. By definition of $\tau_{n}$ there exists a sequence
$\{(s_{j},I_{k_{j}}^{(n)})\dvtx j\geq1\}$ such that
%
%
\begin{equation}
s_{j}\downarrow\tau_{n} \qquad\mbox{as }
j\uparrow\infty\quad\mbox{and}\quad
X_{s_{j}}\bigl(I_{k_{j}}^{(n)}\bigr)\geq2^{-n(1-2\rho)},\qquad j\geq1.
\end{equation}
There exists a subsequence $\{j_{r}\dvtx r\geq1\}$ such that
$I_{k_{j_{r}}}%
^{(n)}=I_{\tilde{k}}^{(n)}$ for some $\tilde{k}\in\mathsf{Z}$. Hence,
for the
fixed $\omega\in A_{n}^{\rho}$,
%
%
\begin{equation}\label{X_tau}%
X_{\tau_{n}}\bigl(I_{\tilde{k}}^{(n)}\bigr) = \lim_{r\rightarrow\infty
}X_{s_{j_{r}}%
}\bigl(I_{\tilde{k}}^{(n)}\bigr)\geq2^{-n(1-2\rho)}.
\end{equation}
Put $\tilde{B}:=[\tilde{k}2^{-n}-2^{-n(1-\rho)},(\tilde{k}+1)2^{-n}%
+2^{-n(1-\rho)}]$. Then there is a constant
$c_{\mazinti{(\ref{int_bound})}}$ independent of $ \rho$ such that
%
%
\begin{eqnarray}\label{int_bound}
\int_{\tilde{B}}dy\, p_{t-s}^{\alpha}(y-z) \geq c_{\mazinti{(\ref{int_bound}%
)}} \hspace*{30pt}\nonumber\\[-8pt]\\[-8pt]
\eqntext{\mbox{for all } z\in I_{\tilde{k}}^{(n)}\mbox{ and } s\in\bigl[
t-2^{-\alpha(1-\rho)n},t\bigr).}
\end{eqnarray}
Now, by the strong Markov property,%
\begin{eqnarray*}
\mathbf{E}X_{t}(\tilde{B})&=& \mathbf{E} e^{a(t-\tau_{n})}%
S_{t-\tau_{n}}^{\alpha}X_{\tau_{n }}(\tilde{B}) \\
&\geq& e
^{-|a|t} \mathbf{E} \biggl\{\int_{\tilde{B}} dy \int_{\mathsf{R}%
} X_{\tau_{n}}(dz)p_{t-\tau_{n}}^{\alpha}(y-z);A_{n}^{\rho} \biggr\}\\
&\geq& e^{-|a|t} \mathbf{E} \biggl\{\int_{I_{\tilde{k}}^{(n)}}%
X_{\tau_{n}}(dz)\int_{\tilde{B}}dy\, p_{t-\tau_{n}}^{\alpha
}(y-z); A_{n}^{\rho} \biggr\}\\
&\geq& c_{\mazinti{(\ref{int_bound})}} \mathbf{E}%
\bigl\{X_{\tau_{n}}\bigl(I_{\tilde{k}}^{(n)}\bigr);A_{n}^{\rho}\bigr\}.
\end{eqnarray*}
Taking into account (\ref{X_tau}) and (\ref{int_bound}) then gives
%
%
\begin{equation}\label{E_bound}%
\mathbf{E}X_{t}(\tilde{B}) \geq c_{\mazinti{(\ref{int_bound})}} 2^{-n(1-2\rho)}
\mathbf{P}(A_{n}^{\rho}).
\end{equation}
On the other hand, in view of Corollary \ref{C1},
%
%
\begin{eqnarray}\label{fin}
\mathbf{E}X_{t}(\tilde{B}) &\leq& |\tilde{B}| \mathbf{E}\sup_{0\leq
x\leq1}X_{t}(x)\nonumber\\[-8pt]\\[-8pt]
&\leq& 2\bigl(2^{-n}+2^{-n(1-\rho)}\bigr) \mathbf{E}\sup_{0\leq x\leq
1}X_{t}(x)\leq C 2^{-n(1-\rho)},\nonumber
\end{eqnarray}
where we wrote $|\tilde{B}|$ for the length of the interval $\tilde{B}$.
Combining (\ref{E_bound}) and (\ref{fin}) completes the proof.
\end{pf}
\begin{lemma}[{[Second term in (\ref{3terms})]}]\label{L.2t}
For fixed $\varepsilon\in(0,\frac{1}{1+\beta})$ and all $n$ large enough,
%
%
\begin{equation}
P_{n,k}^{\varepsilon} \leq2^{-3n/2},\qquad 2\leq k\leq2^{n}+1.
\end{equation}
\end{lemma}
\begin{pf}
Since $T_{-}\leq t_{n}$ on $A^{\varepsilon}$ [recall notation
(\ref{not.tn})],
%
%
\begin{equation}
P_{n,k}^{\varepsilon} \leq\mathbf{P} \Bigl( \sup_{v\leq t_{n}}%
L_{v}\mathrm{1} \Bigl\{\sup_{u\leq v}\Delta L_{u}\leq2^{-n(\eta_{\mathrm{c}%
}+\varepsilon)} \Bigr\}\geq2^{-n\eta_{\mathrm{c}}} \Bigr).
\end{equation}
Applying now Lemma \ref{L3}, with notation of $t_{n}$ from (\ref{not.tn}) we
obtain
%
%
\begin{equation}
P_{n,k}^{\varepsilon} \leq\bigl(c_{\varepsilon}2^{\varepsilon\beta
n-(1-\eta_{\mathrm{c}})(1+\beta)n}+c_{\varepsilon}2^{\eta_{\mathrm{c}}%
(1+\beta)n+\varepsilon\beta n-(\alpha-\beta-\varepsilon)n}%
\bigr)^{ (2^{n\varepsilon})}.
\end{equation}
Inserting the definition of $\eta_{\mathrm{c}}$ and making $n$ sufficiently
large, the estimate in the lemma follows.
\end{pf}

In order to deal with the third term $P_{n,k }^{\varepsilon,\varrho}%
$, we need to define additional events
%
%
\begin{eqnarray}\quad
A_{n,k}^{\varepsilon,\rho,1} &:=& \biggl\{
\mbox{There exists a jump of }M\mbox{ of the form }r^{\ast
}\delta_{(s^{\ast},y^{\ast})}\nonumber\\
&&\hspace*{6.2pt} \mbox{for some }(r^{\ast},s^{\ast},y^{\ast}%
) \mbox{ such that }r^{\ast}\geq(t-s)^{{1}/({1+\beta})+2\varepsilon
/\alpha},\\
&&\hspace*{34.5pt} \biggl\vert\frac
{k+1}{2^{n}}-y^{\ast} \biggr\vert\leq(t-s)^{1/\alpha-2\varepsilon}%
, s^{\ast}\geq t-2^{-\alpha(1+\rho)n} \biggr\} \nonumber
\end{eqnarray}
and
\begin{eqnarray*}
A_{n,k}^{\varepsilon,\rho,2} &:=& A_{n}^{\rho,\mathrm{c}}\cap A_{n,k}%
\\
&&{}\cap\biggl\{\mbox{There exists a jump of }M \mbox{ of the form }%
r^{\ast}\delta_{(s^{\ast},y^{\ast})}\\
&& \hspace*{15.6pt}\mbox{ for some }(r^{\ast},s^{\ast},y^{\ast})\mbox{ such that }%
r^{\ast}\geq(t-s)^{{1}/({1+\beta})+2\varepsilon/\alpha},\\
&&\hspace*{17.6pt} y^{\ast}\in\biggl[\frac{k+1/2}{2^{n}}, \frac{k+1+2^{\rho n+\alpha
2\varepsilon(1-\rho)n}}{2^{n}}\biggr],\\
&&\hspace*{107.65pt} s^{\ast}\in\bigl[t-2^{-\alpha
(1-\rho)n}, t-2^{-\alpha(1+\rho)n}\bigr] \biggr\}.
\end{eqnarray*}
So far we assumed that $ \varepsilon\in(0,\frac{1}{1+\beta})$
and $ \rho\in(0,\frac{1}{2})$. Suppose additionally
that
%
%
\begin{equation} \label{equt:9}%
\frac{\alpha(\alpha+1)2\varepsilon}{1-\eta_{\mathrm{c}}+2\varepsilon
(\alpha^{2}+\alpha-1)} \leq\rho.
\end{equation}
\begin{lemma}[{[Splitting of the third term in (\ref{3terms})]}]\label{lem:9}
For $\rho,\varepsilon>0$ sufficiently small and satisfying (\ref{equt:9})
we have
%
%
\begin{equation} \label{splits}%
P_{n,k}^{\varepsilon,\varrho} \leq\mathbf{P}(A_{n,k}^{\varepsilon,\rho
,1})+\mathbf{P}(A_{n,k}^{\varepsilon,\rho,2})
\end{equation}
for all $ 0\leq k\leq2^{n}-1$ and $n\geq n_{\varepsilon}$.
\end{lemma}
\begin{pf}
First let us describe the strategy of the proof. We are going to show that
whenever a jump of $L_{n,k}^{-}(s), s\leq T_{- }$, of size at least
$2^{-n(\eta_{\mathrm{c}}+\varepsilon)}$ occurs, then it may happen only
in the
points indicated in the definition of $A_{n,k}^{\varepsilon,\rho,1}$ and
$A_{n,k}^{\varepsilon,\rho,2}$. To show this we will in fact show that outside
the sets mentioned in $A_{n,k}^{\varepsilon,\rho,1}$ and $A_{n,k}%
^{\varepsilon,\rho,2}$ the jumps of $L_{n,k}^{-}(s), s\leq T_{- }$, are less
than $2^{-n(\eta_{\mathrm{c}}+\varepsilon)}$.

To implement this strategy, first let us recall that all the jumps of
$ L_{n,k}^{-}(s), s\leq T_{- }$, equal to
%
%
\begin{equation} \label{5.3}%
\Delta X_{s\ast}(y^{\ast}) \biggl(p_{t-s}^{\alpha
} \biggl( \frac{k+1}{2^{n}}-y^{\ast} \biggr) -p_{t-s}^{\alpha} \biggl(
\frac{k}{2^{n}}-y^{\ast} \biggr) \biggr) _{+}
\end{equation}
for some $(s^{\ast},y^{\ast})\in[0,t)\times\mathsf{R}$.

Recall that by definition (\ref{A.eps}), on the event $A^{\varepsilon}$,
%
%
\begin{equation} \label{equt:5}%
\vert\Delta X_{s} \vert\leq c_{\mazinti{(\ref{inL8})}}(t-s)^{(1+\beta
)^{-1}-\gamma}
\end{equation}
with $ \gamma\in(0,(1+\beta)^{-1})$. On the other hand using
Lemma \ref{L1} with $\delta=1$ we obtain
%
%
\begin{equation}\label{equt:6}%
p_{t-s}^{\alpha} \biggl( \frac{k+1}{2^{n}}-y \biggr) -p_{t-s}^{\alpha
} \biggl( \frac{k}{2^{n}}-y \biggr) \leq C2^{-n}(t-s)^{-2/\alpha} .
\end{equation}
From (\ref{equt:5}) and (\ref{equt:6}) we infer
%
%
\begin{eqnarray} \label{equt:7}%
&& \sup_{s\leq t-2^{-\alpha(1-\rho)n}}\Delta X_{s}\sup_{y\in\mathsf{R}} \biggl(
p_{t-s}^{\alpha} \biggl( \frac{k+1}{2^{n}}-y \biggr)
-p_{t-s}^{\alpha} \biggl( \frac{k}{2^{n}}-y \biggr) \biggr)
\nonumber\\
&&\qquad \leq
Cc_{\mazinti{(\ref{inL8})}}2^{-n}\bigl(2^{-\alpha(1-\rho)n}\bigr)
^{{1}/({1+\beta})-\gamma-2/\alpha}
\\
&&\qquad= C 2^{-n(\eta_{\mathrm{c}}-\alpha
\gamma+\rho(1-\eta_{\mathrm{c}}+\alpha\gamma))}.\nonumber
\end{eqnarray}
Furthermore if the jump $\Delta X_{s}$ occurs at the point $y^{\ast}$ with
%
%
\begin{equation}
\biggl\vert y^{\ast}-\frac{k+1}{2^{n}} \biggr\vert\geq(t-s)^{1/\alpha
-2\varepsilon},
\end{equation}
then again by Lemma \ref{L1}, for any $\delta\in[0,1]$,
%
%
\begin{eqnarray}
&&p_{t-s}^{\alpha} \biggl( \frac{k+1}{2^{n}}-y \biggr) -p_{t-s}^{\alpha
} \biggl( \frac{k}{2^{n}}-y \biggr) \nonumber\\[-8pt]\\[-8pt]
&&\qquad\leq C 2^{-n\delta}(t-s)^{-\delta
/\alpha} p_{t-s}^{\alpha}
\bigl((t-s)^{1/\alpha-2\varepsilon}\bigr).\nonumber
\end{eqnarray}
Since
%
%
\begin{equation}
p_{1}^{\alpha}(x) \leq C x^{-1-\alpha},\qquad x\in\mathsf{R},
\end{equation}
we get the bound
%
%
\begin{equation}\hspace*{25pt}
p_{t-s}^{\alpha} \biggl( \frac{k+1}{2^{n}}-y \biggr) -p_{t-s}^{\alpha
} \biggl( \frac{k}{2^{n}}-y \biggr) \leq C 2^{-n\delta}(t-s)^{-
({\delta+1})/{\alpha}+2\varepsilon(\alpha+1)} .
\end{equation}
Hence
%
%
\begin{eqnarray}\hspace*{32pt}
&&\sup_{s<t} \sup_{y\dvtx\vert y-({k+1})/{2^{n}} \vert\geq
(t-s)^{1/\alpha-2\varepsilon}}\Delta X_{s}(y)
\biggl(p_{t-s}^{\alpha} \biggl( \frac{k+1}{2^{n}}-y \biggr)\nonumber\\
&&\hspace*{169.1pt}{} -p_{t-s}^{\alpha
} \biggl( \frac{k}{2^{n}}-y \biggr) \biggr) \\
&&\qquad\leq Cc_{\mazinti{(\ref{inL8})}} 2^{-n\delta}(t-s)^{-({\delta+1})/{\alpha
}+2\varepsilon(\alpha+1)+{1}/({\beta+1})-\gamma}.\nonumber
\end{eqnarray}
Set
%
%
\begin{equation}
\delta:= \eta_{\mathrm{c}}+\alpha\bigl(2\varepsilon(\alpha+1)-\gamma\bigr).
\end{equation}
Note that for all $\varepsilon$ and $\gamma$ sufficiently small, we have
$\delta\in[0,1]$, and we can apply the previous estimates. Thus we
obtain
%
%
\begin{eqnarray} \label{equt:8'}\hspace*{32pt}
&&\sup_{s<t} \sup_{y\dvtx\vert y-({k+1})/{2^{n}} \vert\geq
(t-s)^{1/\alpha-2\varepsilon}}\Delta X_{s}(y)
\biggl(p_{t-s}^{\alpha} \biggl( \frac{k+1}{2^{n}}-y \biggr)\nonumber\\
&&\hspace*{169.1pt}{} -p_{t-s}^{\alpha
} \biggl( \frac{k}{2^{n}}-y \biggr) \biggr)\\
&&\qquad\leq Cc_{\mazinti{(\ref{inL8})}} 2^{-n(\eta_{\mathrm{c}}+\alpha(2\varepsilon
(\alpha+1)-\gamma))}.\nonumber
\end{eqnarray}
Now if we take $ \gamma=2\varepsilon(\alpha+1-1/\alpha)$, which
belongs to these admissible $\gamma$, and $\rho$ as in (\ref{equt:9}), we
conclude that the right-hand side of (\ref{equt:7}) and (\ref{equt:8'}) is
bounded by
%
%
\begin{equation} \label{equt:10'}%
C 2^{-n(\eta_{\mathrm{c}}+2\varepsilon)}.
\end{equation}

For any jump $r^{\ast}\delta_{(s^{\ast},y^{\ast})}$ of $M$ such that
$r^{\ast
}\leq(t-s)^{{1}/({1+\beta})+2\varepsilon/\alpha}$ and $s^{\ast}<t$ we may
apply Lemma \ref{L1} with $\delta=\eta_{\mathrm{c}}+2\varepsilon$ to
get that
%
%
\begin{equation} \label{equt:11'}\quad
\Delta X_{s\ast}(y^{\ast}) \biggl(p_{t-s}^{\alpha
} \biggl( \frac{k+1}{2^{n}}-y^{\ast} \biggr) -p_{t-s}^{\alpha}
\biggl(\frac{k}{2^{n}}-y^{\ast} \biggr) \biggr) \leq C 2^{-n(\eta_{\mathrm{c}%
}+2\varepsilon)}.
\end{equation}
Now recall (\ref{5.3}). Hence combining (\ref{equt:7}), (\ref{equt:8'}),
(\ref{equt:10'}) and (\ref{equt:11'}) the conclusion of Lemma \ref
{lem:9} follows.
\end{pf}

In the next two lemmas we will bound the two probabilities on the right-hand
side of (\ref{splits}).
\begin{lemma}[{[First term in (\ref{splits})]}]\label{L.2split}
For all $ \rho,\varepsilon>0$ sufficiently small and satisfying
%
%
\begin{equation}\label{and}
6 \varepsilon(\alpha+1+\beta) \leq\rho,%
\end{equation}
we have
%
%
\begin{equation}
\mathbf{P}(A_{n,k}^{\varepsilon,\rho,1}) \leq2^{-n-n\rho/2}%
\end{equation}
for all $ k,n$ considered.
\end{lemma}
\begin{pf}
It is easy to see that
\begin{eqnarray*}
A_{n,k}^{\varepsilon,\rho,1} &\subseteq& \bigcup_{l=(1+\rho)n}^{\infty
} \biggl\{\mbox{There exists a jump of }M\mbox{ of the form }r^{\ast}%
\delta_{(s^{\ast},y^{\ast})}\\
&&\hspace*{38.2pt} \mbox{ for some }(r^{\ast},s^{\ast},y^{\ast
})\mbox{ such that }r^{\ast}\geq2^{-l({\alpha}/({1+\beta
})+2\varepsilon)},\\
&&\hspace*{40.2pt} \biggl\vert
\frac{k+1}{2^{n}}-y^{\ast} \biggr\vert\leq2^{-l(1-2\varepsilon\alpha
)}, s^{\ast}\in\bigl[t-2^{-\alpha l}, t-2^{-\alpha(l+1)}\bigr) \biggr\} \\
&=&\!:
\bigcup_{l=(1+\rho)n}^{\infty}A_{n,k,l }^{\varepsilon
,\rho,1}.
\end{eqnarray*}
Recall the random measure $N$ describing the jumps of $ X$. Write
$ Y_{n,k,l}$ for the $N$-measure of
\begin{eqnarray*}
&&\bigl[t(1-2^{-\alpha l}),t\bigl(1-2^{-\alpha(l+1)}\bigr)\bigr]\times\biggl[\frac
{k+1}{2^{n}}-2^{-l(1-2\alpha\varepsilon)}, \frac
{k+1}{2^{n}}+2^{-l(1-2\alpha
\varepsilon)}\biggr]\\
&&\qquad{}\times\bigl[2^{-l({\alpha}/({1+\beta})+2\varepsilon
)}, \infty\bigr).
\end{eqnarray*}
Then, by Markov's inequality,
%
%
\begin{equation}
\mathbf{P}(A_{n,k,l}^{\varepsilon,\rho,1}) = \mathbf{P}(Y_{n,k,l}%
\geq1) \leq\mathbf{E}Y_{n,k,l }.
\end{equation}
Therefore,
%
%
\begin{equation}
\mathbf{P}(A_{n,k}^{\varepsilon,\rho,1}) \leq\sum_{l\geq(1+\rho
)n}\mathbf{P}(A_{n,k,l}^{\varepsilon,\rho,1}) \leq\sum_{l\geq(1+\rho
)n}\mathbf{E}Y_{n,k,l }.
\end{equation}
From the formula for the compensator of $N$ we get
%
%
\begin{eqnarray}
\mathbf{E}Y_{n,k,l} &=& \varrho\int_{t(1-2^{-\alpha l})}^{t(1-2^{-\alpha
(l+1)})}ds\, \mathbf{E}X_{s} \biggl( \biggl[\frac{k+1}{2^{n}%
}-2^{-l(1-2\alpha\varepsilon)},\nonumber\\
&&\hspace*{118pt}\frac{k+1}{2^{n}}+2^{-l(1-2\alpha
\varepsilon
)} \biggr] \biggr)\nonumber\\[-8pt]\\[-8pt]
&&{} \times\int_{2^{-l({\alpha}/({1+\beta})+2\varepsilon)}%
}^{\infty}dr\, r^{-2-\beta}\nonumber\\
&\leq& C 2^{-\alpha l} 2^{-l(1-2\alpha\varepsilon)} 2^{l(\alpha
+2\varepsilon(1+\beta))}.\nonumber
\end{eqnarray}
Consequently,
%
%
\begin{equation}\qquad
\mathbf{P}(A_{n,k,l}^{\varepsilon,\rho,1}) \leq C\sum_{l\geq(1+\rho
)n}2^{-l+2\varepsilon(\alpha+1+\beta)l} \leq C 2^{-(1+\rho
)n+2\varepsilon
(\alpha+1+\beta)(1+\rho)n}.
\end{equation}
Noting that $2\varepsilon(\alpha+1+\beta)(1+\rho)\leq\rho/2$ under the
conditions in the lemma, we complete the proof.
\end{pf}
\begin{lemma}[{[Second term in (\ref{splits})]}]\label{lem:10}
For all $\varepsilon,\rho>0$ sufficiently small,
%
%
\begin{equation}
\mathbf{P}(A_{n,k}^{\varepsilon,\rho,2}) \leq2^{-3n/2}%
\end{equation}
for all $ k,n$ considered.
\end{lemma}
\begin{pf}
It is easy to see by construction that
\begin{subequation}
\label{big.array}%
\begin{eqnarray}
&&A_{n,k}^{\varepsilon,\rho,2}\subseteq A_{n}^{\rho,\mathrm{c}}\cap \biggl\{
\mbox{There exist at least two jumps of }M\nonumber\\
&&\hspace*{79.1pt} \mbox{of the form }r^{\ast}\delta_{(s^{\ast},y^{\ast}%
)}\mbox{ such that}\nonumber\\
\label{a}
&&\hspace*{79.1pt} r^{\ast}\geq2^{-n({\alpha(1+\rho)}/({1+\beta})+2\varepsilon
(1+\rho))},\\
\label{b}
&&\hspace*{79.1pt} y^{\ast}\in\biggl[\frac{k-2}{2^{n}}, \frac{k+1+2^{\rho
n+2\alpha\varepsilon(1-\rho)n}}{2^{n}} \biggr] ,\\
\label{c}%
&&\hspace*{79.1pt}\hspace*{11.46pt} s^{\ast}\in\bigl[ t-2^{-\alpha(1-\rho)n}, t-2^{-\alpha(1+\rho
)n} \bigr] \biggr\}.
\end{eqnarray}
On the event $A_{n}^{\rho,\mathrm{c}}$, for the intensity of jumps satisfying
(\ref{a})--(\ref{c}), we have
\end{subequation}
\begin{eqnarray*}
&& \int_{t-2^{-\alpha(1-\rho)n}}^{t-2^{-\alpha(1+\rho)n}} d%
s\, X_{s} \biggl( \biggl[\frac{k-2}{2^{n}}, \frac{k+1+2^{\rho n+2\alpha
\varepsilon(1-\rho)n}}{2^{n}}\biggr] \biggr)\\
&&\quad{}\times \int_{2^{-n(
{\alpha(1+\rho)}/({1+\beta})+2\varepsilon(1+\rho))}}^{\infty}d%
r\, r^{-2-\beta}\\
&&\qquad \leq2^{-\alpha(1-\rho)n} 2^{-n(1-2\rho)} 2^{\rho n+2\alpha
\varepsilon(1-\rho)n+2} 2^{n(\alpha(1+\rho)+2\varepsilon(1+\rho
)(1+\beta))}\\
&&\qquad \leq2^{-n} 2^{10(\rho+2\varepsilon)n} \leq2^{-{3}/{4}n}%
\end{eqnarray*}
for all $\varepsilon$ and $\rho$ sufficiently small. Since the number
of such
jumps can be represented by means of a time-changed standard Poisson process,
the probability to have at least two jumps is bounded by the square of the
above bound and we are done.
\end{pf}
\begin{lemma}[{[Third term in (\ref{3terms})]}]\label{L.3t}For all $\rho
,\varepsilon>0$ sufficiently small, satisfying
(\ref{equt:9}) and (\ref{and}), we have
%
%
\begin{equation}
P_{n,k}^{\varepsilon,\varrho} \leq2^{-3n/2}+C2^{-n-\rho n/2},\qquad 2\leq
k\leq2^{n}+1, n\geq n_{\varepsilon}.
\end{equation}
\end{lemma}
\begin{pf}
The proof follows immediately from Lemmas \ref{lem:9}, \ref{L.2split}
and \ref{lem:10}.
\end{pf}
\begin{pf*}{Proof of Lemma \ref{prop:2}}
Applying Lemmas \ref{L.small.prob}, \ref{L.2t} and \ref{L.3t} to
(\ref{3terms}) we obtain
%
%
\begin{equation}\hspace*{28pt}
\mathbf{P} \Biggl( \bigcup_{k=2}^{2^{n}+1} ( A^{\varepsilon}\cap
A_{n,k}\cap B_{n,k}^{-,\mathrm{c}} ) \Biggr) \leq
c_{\mazinti{(\ref{small.prob})}} 2^{-\rho n}+2^{-n/2}+C2^{-\rho n/2}+2^{-n/2}%
\end{equation}
for all $\rho,\varepsilon>0$ sufficiently small satisfying (\ref{equt:9})
and (\ref{and}) as well as all $n\geq n_{\varepsilon}$. Since these terms
are summable in $n$, the claim of the lemma follows.
\end{pf*}
\begin{pf*}{Proof of Lemma \ref{Prop2}}
The proof follows immediately from (\ref{5.10}) and
Lemmas \ref{L.1part} and \ref{prop:2}.
\end{pf*}

\section*{Acknowledgments}
We thank Don Dawson for very helpful discussions of the subject. We
also thank the referees for their useful comments and suggestions which
improved the exposition.

%

%
\printaddresses

\end{document}